\documentclass[leqno,12pt]{article}
\usepackage{amssymb,amsmath,graphicx}

\title{Local Euler-Maclaurin formula for polytopes}
\author{Nicole Berline and Mich{\`e}le Vergne}
\date{\today}

\begin{document}
\maketitle
\newtheorem{theorem}{Theorem}
\newtheorem{proposition}[theorem]{Proposition}
\newtheorem{lemma}[theorem]{Lemma}
\newtheorem{definition}[theorem]{Definition}
\newtheorem{corollary}[theorem]{Corollary}
\newtheorem{remark}[theorem]{Remark}
\newtheorem{example}[theorem]{Example}
\newenvironment{proof}{{\bf Proof. }}{\par}
\newcommand{\half}{{\frac{1}{2}}}
\newcommand{\C}{{\mathbb C}}
\newcommand{\R}{{\mathbb R}}
\newcommand{\Z}{{\mathbb Z}}
\newcommand{\N}{{\mathbb N}}
\newcommand{\Q}{{\mathbb Q}}
\newcommand{\CA}{{\cal A}}
\newcommand{\CC}{{\cal C}}
\newcommand{\CE}{{\cal E}}
\newcommand{\CF}{{\cal F}}
\newcommand{\CH}{{\cal H}}
\newcommand{\CL}{{\cal L}}
\newcommand{\CM}{{\cal M}}
\newcommand{\CP}{{\cal P}}
\newcommand{\CS}{{\cal S}}
\newcommand{\CV}{{\cal V}}
\newcommand{\la}{{\langle}}
\newcommand{\ra}{{\rangle}}
\newcommand{\proj}{\operatorname{proj}}
\newcommand{\Res}{\operatorname{Res}}
\newcommand{\Convex}{\operatorname{P}}
\newcommand{\simplify}{\operatorname{Simplify}}
\newcommand{\Sum}{\operatorname{Sum}}
\newcommand{\card}{\operatorname{Card}}
\newcommand{\dir}{\operatorname{dir}}
\newcommand{\linapex}{\operatorname{linapex}}
\newcommand{\lin}{\operatorname{lin}}
\newcommand{\apex}{\operatorname{apex}}
\newcommand{\relvol}{\operatorname{rvol}}
\newcommand{\vertex}{\operatorname{vert}}
\newcommand{\Vertex}{\operatorname{Vertex}}
\newcommand{\cone}{\operatorname{cone}}
\newcommand{\normal}{\operatorname{No}}
\newcommand{\Ev}{\operatorname{Ev}}

\newcommand{\Ber}{\operatorname{Ber}}
\newcommand{\vol}{\operatorname{vol}}
\renewcommand{\a}{{\mathfrak{a}}}
\renewcommand{\c}{{\mathfrak{c}}}
\renewcommand{\d}{{\mathfrak{d}}}
\newcommand{\f}{{\mathfrak{f}}}
\newcommand{\n}{{\mathfrak{n}}}
\renewcommand{\t}{{\mathfrak{t}}}
\newcommand{\co}{{\mathfrak{c}\mathfrak{o}}}

\newcommand{\p}{{\mathfrak{p}}}
\newcommand{\lattice}{\Lambda}
\section{Introduction}
By the name Euler-Maclaurin, one refers to formulas which  relate
discrete sums to integrals, in particular in the following
framework. Let $\p$ be a rational convex polytope in $\R^d$ and let
$h(x)$ be a polynomial  function on $\R^d$. The sum of the values
$h(x)$ at integral points of $\p$ is written as a sum of terms
indexed by the set $\CF(\p)$  of faces of $\p$,
\begin{equation}\label{intro-maclaurin}
\sum_{x\in \p\cap\Z^d}h(x)=\sum_{\f\in\CF(\p)}\int_\f
D(\p,\f)\cdot h
\end{equation}
where,  for  each face $\f$ of $\p$ , $D(\p,\f)$ is a differential
operator (of infinite order) with constant coefficients on $\R^d$.
The basic example is the historical Euler-Maclaurin summation
formula in dimension 1: for $a_1\leq a_2 \in \Z$,
\begin{eqnarray*}
\sum_{a_1}^{ a_2}h(x) = \int_{a_1}^{a_2}h(t)dt -  \sum_{n\geq
1}\frac{b(n)}{n!}h^{(n-1)}(a_1) +\sum_{n\geq 1}(-1)^{n}
\frac{b(n)}{n!} h^{(n-1)}(a_2)
\end{eqnarray*}
where $b(n)$ are the  Bernoulli numbers.

When $\p$ is an integral polytope, the existence of such operators
is the combinatorial counterpart of a homological property of the
associated toric variety: the invariant cycles corresponding to
the faces of $\p$ generate the equivariant homology module. An
Euler-Maclaurin formula amounts to an explicit Riemann-Roch
theorem, as obtained by Khovanskii and Pukhlikov \cite{KP92} for
integral polytopes corresponding to smooth toric varieties, and
extended by Cappell and Shaneson \cite{cs1} to any integral
polytope. Furthermore, by transforming Cappell-Shaneson
homological methods into  purely combinatorial techniques, valid
for any {rational} (not necessary integral) polytope, Brion and
Vergne \cite{bv1} obtained various expressions for the sum
$\sum_{x\in \p\cap\Z^d}h(x)$, either as an integral over a
deformed polytope followed by differentiation with respect to the
deformation parameter, or formulas of type
(\ref{intro-maclaurin}).

In this article, we construct differential operators $D(\p,\f)$,
with rational coefficients, which satisfy  Formula (\ref
{intro-maclaurin}) and which moreover enjoy two essential
properties: they are \emph{local} and they are \emph{computable}.
The existence of operators with these properties was conjectured
in \cite{barpom}.

By {local}, one  means that  $D(\p,\f)$ depends only on the
equivalence class - modulo integral translations - of the
{\emph{transverse cone}  $\t(\p,\f)$ of $\p$ along $\f$ (see
Definition \ref{transverse}  and Figure \ref{normaldim3}). In
particular, if $\p$ is an integral polytope, the operator $D(\p,\f)$
depends only on the cone of  feasible directions of $\p$ along $\f$.

By {computable}, one means that there exists an algorithm which
computes the  $m$ lowest order terms of $D(\p,\f)$, with running
time polynomial with respect to the size of the data defining $\p$,
at least  when the dimension $d$  and the number $m$ are fixed.

On the contrary, Cappell-Shaneson and Brion-Vergne operators are
neither local nor computable.

When applied to  the constant polynomial $h(x)=1$, Formula
(\ref{intro-maclaurin}) takes the form
\begin{equation}\label{intro-danilov}
{\rm Card} (\p\cap\Z^d) =\sum_{\f\in\CF(\p)}\nu_0(\p,\f) \vol(\f)
\end{equation}
where the coefficients $\nu_0(\p,\f)$ are rational numbers. In the
context of toric varieties, Danilov \cite{danilov} asked the
question of existence of coefficients $\nu_0(\p,\f)$ with the local
property, in the case of an integral polytope. This result was
proven by Morelli \cite{morelli} and McMullen \cite{McMullen}. In
\cite{pomm}, Pommersheim and Thomas gave a canonical construction of
rational coefficients $\nu_0(\p,\f)$ which satisfy Formula
(\ref{intro-danilov}), as a consequence of their expression for the
Todd class of a toric variety. In a companion article, we will
similarly obtain a local formula  for the \textbf{equivariant} Todd
class of any toric variety. The result is stated in Theorem
\ref{todd}.

\medskip
The computability of our operators $D(\p,\f)$ extends the
following remarkable result of Barvinok \cite{barvinok}: when the
dimension $d$ is fixed, the number of integral points ${\rm Card}
(\p\cap\Z^d)$
 can be computed by a polynomial time  algorithm.
 This  result  is a consequence of Brion's theorem
\cite{b}, according to which the computation can be distributed
over the tangent cones at the vertices,  and of
 Barvinok's signed decomposition of a cone into unimodular cones by a polynomial time algorithm.
Based on this method, efficient software packages for  integer
points counting problems and the effective computation of $\sum_{x
\in \p \cap \Z^d} h(x)$ have been developed \cite{latte},
 \cite{deloera}, \cite{ver-woo}, \cite{ver-seg}.
\medskip

Let us now explain the construction of the operators $D(\p,\f)$. We
define the differential operator $D(\p,\f)$ through its symbol,
using a scalar product on $\R^d$. For this purpose, we first
associate an analytic function $\mu(\a)$ (defined on a neighborhood
of 0) on $\R^d$ to any rational affine cone $\a\subset \R^d$.

For instance, for a half-line $s+\R_+$, we have
$$
\mu(s+\R_+)(\xi)=\frac{e^{[[s]]\xi}}{1-e^\xi}+\frac{1}{\xi}.
$$
where  $[[s]]=n-s$, with $n\in\Z$ , $n-1<s\leq n$. Remark that
this function is analytic at $0$, with value at $0$ given by
$\mu(s+\R_+)(0)=\frac{1}{2}-[[s]]$.

The assignment $\a\mapsto \mu(\a)$ has beautiful geometric
properties. The most important one is that it is  a
\textbf{valuation} when the vertex of $\a$ is fixed (Theorem
\ref{maintheoremplus}). For instance,
\begin{equation}\label{val}
\mu(\a_1\cup \a_2)=\mu(\a_1)+ \mu(\a_2)-\mu(\a_1\cap \a_2).
\end{equation}

Moreover, $\mu(\a)$ is unchanged when $\a$ is moved by a lattice
translation, and  the map $\a\mapsto \mu(\a)$ is equivariant with
respect to lattice-preserving isometries. In particular, with the
standard scalar product,  we thus get invariants of the group
${\rm O}(n,\Z)$.

\medskip
If $\a$ contains a straight line, we set $\mu(\a)=0$. If $\a$ is
 pointed  with vertex $s$, we define $\mu(\a)$ recursively by the
relation:
\begin{equation}\label{intro-mu}
\mu(\a)(\xi)= e^{-\langle \xi,s\rangle }\left(\sum_{x\in
\a\cap\Z^d}e^{\langle \xi,x\rangle }\,+ \, \sum_{\f,
\dim(\f)>0}\mu(\t(\a,\f))(\xi)\int_\f e^{\langle \xi,x\rangle }
dm_\f(x)\right)
\end{equation}
where $\f$ denotes a face of $\a$ and $dm_\f(x)$ denotes the
canonical Lebesgue measure on $\f$ defined by the lattice. In
(\ref{intro-mu}),  the function $\mu(\t(\a,\f))$ is a priori defined
only on a subspace of $\R^d$, namely the orthogonal to the face
$\f$.  We extend it to  $\R^d$ by orthogonal projection. Our main
point is to show that Formula (\ref{intro-mu}) actually defines an
analytic function. Thanks to the {valuation property}, the proof is
reduced to the case of a simplicial unimodular cone.

\medskip
Then we define $D(\p,\f)$  as the differential operator with symbol
$\mu(\t(\p,\f))(\xi)$. If the scalar product is rational, the Taylor
series of  $\mu(\a)$ has rational coefficients, in particular the
numbers $\nu_0(\p,\f)$ in (\ref{intro-danilov}) are rational. Note
that $D(\p,\f)$ involves only differentiation in directions
perpendicular  to $\f$.

With this definition, Euler-Maclaurin formula
(\ref{intro-maclaurin}) for any rational polytope $\p$ follows
easily from Brion's theorem. Indeed, the defining formula
(\ref{intro-mu}) is formally Formula (\ref{intro-maclaurin}) where
the polytope $\p$ is replaced by the cone $\a$ and the polynomial
$h(x)$ is replaced by the exponential $e^{\langle \xi,x\rangle }$.

The computability of the functions $\mu(\t(\p,\f))$ is also deduced
from Barvinok's fast decomposition of cones, thanks to the valuation
property.

 Moreover, Barvinok proved recently \cite{bb} that, given an integer $m$, there
exists a polynomial time algorithm which computes the $m$ highest
coefficients of the Ehrhart quasipolynomial of any rational
simplex in $\R^d$, \emph{when the dimension $d$ is considered as
an input}. We hope that our construction of the functions
$\mu(\a)$ will lead to another polynomial time algorithm which
would compute the $m$ highest coefficients of the Ehrhart
quasipolynomial for any simplex in $\R^d$ and any polynomial
$h(x)$, when the dimension $d$ and the degree of $h$ are
considered as input. We want to point out that our construction
involves only cones of dimension less than $m$ when computing the
$m$ highest order Ehrhart coefficients.

In the forthcoming article \cite{bal-ber-ver}, we will compare our
construction to the mixed valuation method of \cite{bb}.

\medskip
Let us illustrate the results in dimension 2. For an affine cone
$\a$ with integral vertex $s$ and edges generated by two integral
vectors $v_1,v_2$ with $\det(v_1,v_2)=1$, (that is to say, $\a$ is
unimodular), we have:
$$
\mu(\a)(0)=
\frac{1}{4}+\frac{\langle v_1,
v_2\rangle}{12}(\frac{1}{\langle v_1,v_1\rangle}
+\frac{1}{\langle v_2,v_2\rangle}).
$$
For a  general cone, we compute $\mu(\a)(0)$ using the valuation
property (\ref{val}).
\begin{figure}[!h]
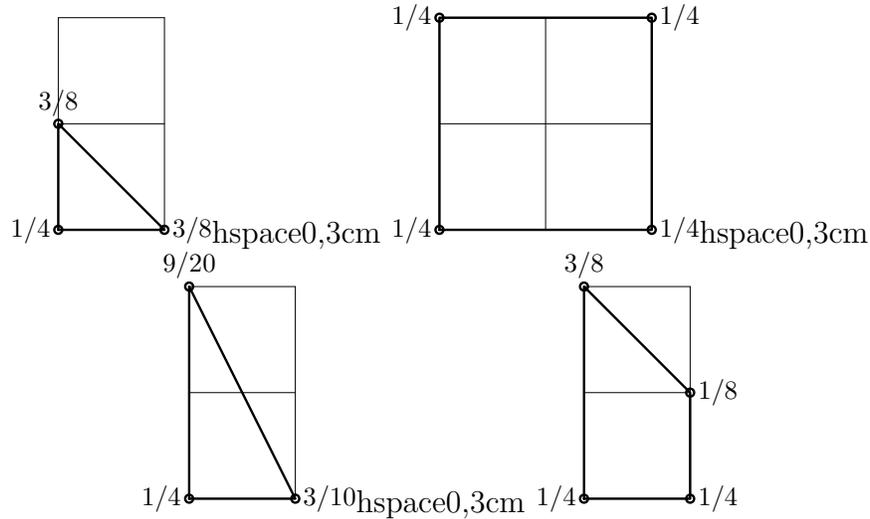

\begin{center}
\includegraphics{delzant-Nicole.01.eps}hspace{0,3cm}
\includegraphics{delzant-Nicole.06.eps}hspace{0,3cm}
\includegraphics{delzant-Nicole.04.eps}hspace{0,3cm}
\includegraphics{delzant-Nicole.02.eps}
\caption{Coefficient $\nu_0(\p,s)=\mu(\t(\p,s))(0)$ at the vertex
$s$}\label{intro-delzant}
\end{center}
\end{figure}
\begin{figure}[!h]
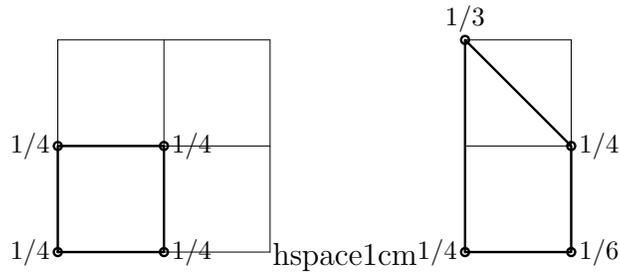

\begin{center}
\includegraphics{delzant-Nicole.07.eps}hspace{1cm}\includegraphics{delzant-Nicole.05.eps}
\caption{ Cappell-Shaneson coefficients}\label{intro-Cappell}
\end{center}
\end{figure}

According to Pick's theorem, the number of integral points of an integral polygon $\p\subset \R^2$ is given by
\begin{equation*}\label{intro-pick}
{\rm Card} (\p\cap\Z^2)={\rm area(\p)} + \frac{1}{2}{\rm
length}_{\Z^2} \,\partial(\p)\,+ \,1.
\end{equation*}
According to our local formula, as well as Pommersheim-Thomas's, we have:
$$
{\rm Card} (\p\cap\Z^2)={\rm area(\p)} + \frac{1}{2}{\rm
length}_{\Z^2} \,\partial(\p)\,+ \sum_s \nu_0(\p,s)
$$
with $\nu_0(\p,s)=\mu(\t(\p,s))(0)$, where $s$ runs over the
vertices of $\p$. Thus the constant $1$ is canonically distributed
over the vertices (Figure \ref{intro-delzant}). Cappell-Shaneson
coefficients \cite{cs2} give a different, non local, distribution of
the constant $1$ over the vertices. For instance, the bottom right
coefficients in the square and the trapezoid  of Figure
\ref{intro-Cappell} are different, although the vertices have the
same tangent cone.

More examples are given at the end of the paper, for polygons $\p$
with rational (non integral) non unimodular  vertices. Based on our
Euler-Maclaurin operators, we wrote a Maple program which computes
the value of the sum  $\sum_{x\in \p\cap\Z^2}x_1^{m_1}x_2^{m_2}$ and
also the (periodic) coefficients of the corresponding Ehrhart
quasipolynomial.

\section{Definitions and notations}

We consider  a rational vector space $V$, that is to say  a finite
dimensional real vector space with a lattice denoted by $\lattice_V$
or simply $\lattice$. By  lattice, we  mean a discrete additive
subgroup of $V$ which  generates  $V$ as a vector space. Hence, a
lattice is  generated by a basis of the vector space $V$. A basis of
$V$ which is a $\Z$-basis of $\lattice_V$ is called an integral
basis.

 We will need to
consider subspaces and quotient spaces of $V$, this is why we
cannot just let $V=\R ^d$ and $\lattice = \Z^d$. The points of
$\lattice$ are called integral. A point $x\in V$ is called
rational if $qx\in \lattice $ for some integer $q\neq 0$. The
space of rational points in $V$ is denoted by $V_\Q$. A subspace
$W$ of $V$ is called rational if $W\cap \lattice $ is a lattice in
$W$.  If $W$ is a  rational subspace, the image of $\lattice$ in
$V/W$ is a lattice in $V/W$, so that $V/W$ is a rational vector
space.

A rational space $V$, with lattice $\lattice$, has a canonical
Lebesgue measure, for which  $V/\lattice$ has measure $1$. An affine
subspace $W$ of $V$ is called rational if it is  a translate of a
rational subspace by a rational element. It is similarly provided
with a canonical Lebesgue measure. We will sometimes denote this
measure  by $dm_W$. For example, let $W$ be a rational line of the
form $W= s+ \R v$. Assume that  $v$ is a generator of the group $\R
v \cap \lattice$ (we say that $v$ is a \emph{primitive} vector).
Then $dm_W( s + tv)= dt$.

If $v_i\in V_\Q$ are linearly independent vectors, we denote by
 $\Box(v_1,\dots, v_k)$ the semi-open parallelepiped generated by
 the $v_i$'s:
$$
\Box(v_1,\dots, v_k)=\sum_{i=1}^k[0,1[v_i.
$$
We denote by  $\vol(\Box(v_1,\dots, v_k))$ its relative volume, that
is to say its  volume with respect to the  canonical measure on the
subspace generated by $v_1,\dots, v_k$.

We denote by $V^* $ the dual space of $V$.  We will denote elements of
$V$ by latin letters $x,y,v,\dots$
and elements of $V^*$ by greek letters $\xi,\alpha,\dots$.  We denote
the duality bracket by  $\langle\xi,x\rangle$.

$V^* $ is equipped with the dual lattice $\lattice^* $
of $\lattice$ :
$$
\lattice^* = \{\xi\in V^*\; ; \,\,\langle\xi,x\rangle \in \Z \; \mbox{for
  all}\;
  x\in \lattice \}.
$$

 If $S$
is a subset of $V$, we denote by $S^{\perp}$ the subspace of $V^*$
orthogonal to $S$:
$$
S^{\perp}= \{\xi\in V^*\; ;\,\,\langle\xi,x\rangle =0 \;\mbox{for all}\;
  x\in  S\}.
$$

If $W$ is a subspace of $V$, the dual space  $(V/W)^*$ is
canonically identified with the subspace $W^{\perp} \subset V^*$.

If $S$ is a subset of $V$, we denote by $<S>$ the affine subspace
generated by $S$. If $S$ consists of rational points, then $<S>$ is
rational. Remark that $<S>$ may  contain no integral point.
 We denote by $\lin(S)$ the vector subspace of $V$ parallel to
$<S>$.

\begin{figure}[!h]
\begin{center}
\includegraphics{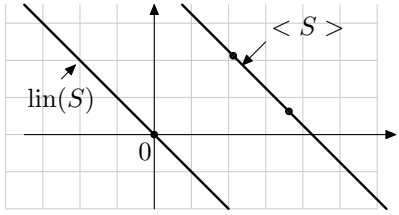}
\caption{The affine space $<S>$ and the linear space
$\lin(S)$}\label{lin}
\end{center}
\end{figure}

The set of non  negative real numbers is denoted by $\R_+$. A  convex rational cone
 $\c$ in $V$ is  a closed
convex cone $\sum_{i=1}^k\R_+ v_i$ which is generated
 by a finite number of elements  $v_i$ of $V_\Q$.
In this article, we simply say  cone instead of convex rational cone.

An affine (rational) cone $\a$ is, by definition, the translate of
a cone in $V$ by an element $s\in
  V_\Q$. This cone is uniquely defined by $\a$; it is called the cone
  of directions of  $\a$ and denoted by $\dir (\a)$. Thus
$\a=s+\dir(\a)$.

 A  cone $\c$ is called simplicial if it is  generated by
independent elements of $V_\Q$. A simplicial cone $\c$ is called unimodular
if it is generated by independent integral vectors  $v_1,\dots, v_k$
such that $\vol(\Box(v_1,\dots, v_k))= 1$.  An  affine cone $\a$ is
called simplicial (resp. simplicial unimodular) if $\dir(\a)$ is
simplicial (resp. simplicial unimodular).

An affine cone $\a$ is called pointed if it does not contain any
straight line.

The set of faces of an affine cone $\a$ is denoted by $\CF(\a)$. If
$\a$ is pointed, then the vertex of $\a$ is the unique face of
dimension $0$, while $\a$ is the unique face of maximal dimension
$\dim \a$.

The dual cone $\c^*$ of a cone $\c$ is the set of $\xi \in V^*$ such
that $\la \xi,x\ra \geq 0$ for any $x\in \c$.

A convex rational polyhedron $\p$ in $V$ (we will simply say
 polyhedron) is, by definition, the intersection of a finite number of
half spaces with boundary a rational affine hyperplane.

\begin{definition}
 We say that
 $\p$ is solid (in V) if $<\p>= V$.
\end{definition}

  The set of faces of
$\p$ is denoted by $\CF(\p)$ and the set of vertices of $\p$ is
denoted by $\CV(\p)$.
\begin{figure}[!h]
\begin{center}
\includegraphics{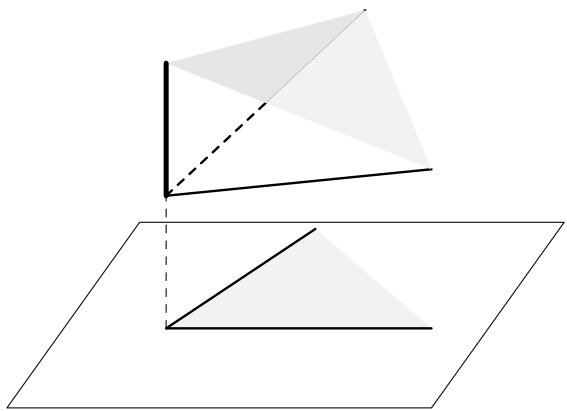}
\caption{The transverse cone along an edge in dimension
$3$}\label{normaldim3}
\end{center}
\end{figure}

We now introduce the main geometrical object in our study, the
\emph{transverse cone of a polyhedron} $\p$ along one of its faces
$\f$ (see Figure \ref{normaldim3}). Let $x$ be a point in the
relative interior of $\f$. Recall that the cone of feasible
directions of $\p$  at $x$   is the set $\c(\p,\f):=\{v\in V\;;
x+\epsilon v\in \p \;\mbox{for}\; \epsilon>0\; \mbox{small enough}
\; \}$. It does not depend on the choice of  $x$ and contains the
linear space $\lin(\f)$. The supporting cone of $\p$ along  $\f$ is
the affine cone $<\f> +\c(\p,\f)$. We denote the projection $V\to
V/\lin(\f)$ by $\pi_\f$.
\begin{definition}\label{transverse}
Let $\p$ be a polyhedron and  $\f$ a face of $\p$. The transverse
cone $\t(\p,\f)$ of $\p$ along $\f$ is   the image $\pi_\f(<\f>
+\c(\p,\f))$ of the supporting cone
 in  $V/\lin(\f)$.
 \end{definition}
We will often write simply  \emph{transverse cone along }$\f$, when
$\p$ is understood.

If $v$ is a vertex of $\p$, the transverse  cone $\t(\p,v)$
coincides with the supporting  cone $v+\c(\p,v)\subset V$.

The transverse cone $\t(\p,\f)$ is a pointed \textbf{affine} cone in
the quotient space $V/\lin(\f)$. Its dimension is equal to the
codimension of $\f$ in $<\p>$.  Its vertex is the projection
$\pi_\f(x)$ of any point $x$ of $\f$ on $V/\lin(\f)$.

If $\a$ is an affine cone and $\f$ is a face of $\a$,  then the
supporting cone is  $ \a + \lin(\f)$ and the transverse cone
$\t(\a,\f)$ along $\f$ is just  the projection $\pi_\f(\a)$  of $\a$
on $V/\lin(\f)$.

\medskip
We shall make use of subdivisions of cones.
\begin{definition}
A subdivision of a cone $\c$ is a finite collection $\cal C$ of
cones in $V$ such that:

(a) The faces of any  cone  in $\cal C$ are in $\cal C$.

 (b) If $\d_1$ and $\d_2$ are two elements of $\cal C$, then the
 intersection $\d_1\cap \d_2$ is a face  of both  $\d_1$ and $\d_2$.

 (c) We have  $\c=\cup_{\d \in \cal C}\d.$

If furthermore the elements of $\cal C$ are simplicial  cones, the
subdivision will be called simplicial.
\end{definition}
\begin{example}\label{subdivisionR}
The basic example is the subdivision
$\{\R_+,(-\R_+), \{0\}\}$ of the one-dimensional cone $\R$.
\end{example}
It is easy to see that any pointed cone admits a subdivision into
simplicial unimodular cones.

\medskip

If $S$ is a subset of $V$, we denote by $\chi(S)$ the characteristic
function of $S$ (also called the indicator function of $S$).

\medskip

Explicit expression of operators  $D(\p,\f)$  involves the Bernoulli
polynomials $b(n,t)$, defined by the generating series
\begin{equation}\label{bernoulli}
\sum_{n=0}^\infty \frac{b(n,t)}{n!} X^n = \frac{e^{t X}X}{e^X-1}.
\end{equation}
The Bernoulli number is $b(n,0)$ and is denoted by $b(n)$.

\section{Meromorphic functions associated to polyhedra}

By a meromorphic function on $V^*$ with rational coefficients, we
mean a meromorphic function on the complexification of $V^*$ which
can be written as the quotient of two holomorphic functions with
rational Taylor coefficients with respect to an integral basis of
$V^*$.

We recall the construction of two meromorphic functions with
rational coefficients on  $V^*$, associated to any polyhedron $\p$
in $V$ (see the survey \cite{barpom}).  The first function $I(\p)$
is defined  via integration over $\p$, the second function $S(\p)$
via summation over the set of  integral points of $\p$.

We denote by $dm_{<\p>}$  the relative Lebesgue measure on the
affine space spanned by $\p$.

\begin{proposition}\label{valuationI}
There exists a map $I$ which to every polyhedron $\p\subset V$
associates a meromorphic function with rational coefficients $I(\p)$
on  $V^*$, so that the following properties hold:

(a) If $\p$ contains a straight line, then $I(\p)$=0.

(b) If $\xi\in V^*$ is such that $|e^{\la \xi,x\ra}|$ is integrable
over $\p$,  then
$$
I(\p)(\xi)= \int_\p e^{\la \xi,x\ra} dm_{<\p>}(x).
$$

(c) For every point $s\in V_\Q$, we have
$$
I(s+\p)(\xi) = e^{\la \xi,s\ra}I(\p)(\xi).
$$

(d) The map $I$ is a {\em solid valuation}: if the characteristic
functions $\chi(\p_i)$ of a family of  polyhedra $\p_i$ satisfy a
linear relation $\sum_i r_i \chi(\p_i)=0$, then the functions
$I(\p_i)$ satisfy the  relation
$$
\sum_{\{i,<\p_i>=V\} }r_i I(\p_i)=0.
$$
\end{proposition}

\begin{example}.

\noindent$\bullet$ If $\p=\{s\}$ is a point, then $I(\p)(\xi)
=e^{\la\xi,s\ra}$.

\noindent$\bullet$ In dimension $1$, if  $\p=s+\R^+$, where $s\in
\Q$, then $I(\p)(\xi)=-\frac{e^{\xi s}}{\xi}.$

\noindent$\bullet$ If $\c$ is a simplicial cone generated by $k\leq
d$ independent vectors $v_1,\dots,v_k$, we have
\begin{equation}\label{Iparallelepiped}
I(\c)(\xi)= (-1)^k  \frac{\vol(\Box(v_1,\dots,v_k)) }{\prod_{i=1}^k
  \la \xi, v_i\ra}.
\end{equation}
\end{example}
These formulas follow immediately from the computation in
dimension 1, $\int_{0}^{\infty}e^{\xi x}dx=\frac{-1}{\xi}$ for
$\xi<0$. Thus,  for a simplicial cone, the function $I(\c)$ can
indeed be extended as a rational function on the whole space
$V^*$.

In a similar way, one defines  the second meromorphic function,
which is the discrete analogue of $I(\p)$.

\begin{proposition}\label{valuationS}
There exists a map $S$ which to every polyhedron $\p\subset V$
associates a meromorphic function with rational coefficients $S(\p)$
on $V^*$ so that the following properties hold:

(a) If $\p$ contains a straight line, then $S(\p)$=0.

(b) If $\xi\in V^*$ is such that $\sum_{x\in \p\cap \lattice}
|e^{\la \xi,x\ra}| <\infty$,  then
\begin{equation}
S(\p)(\xi)= \sum_{x\in \p\cap \lattice} e^{\la \xi,x\ra}.
\end{equation}

(c) For every integral point $s\in \lattice$, we have
$$
S(s+\p)(\xi) = e^{\la \xi,s\ra}S(\p)(\xi).
$$

(d) The map $S$ is a valuation: if the characteristic functions
$\chi(\p_i)$ of a family of  polyhedra $\p_i$ satisfy a linear
relation $\sum_i r_i \chi(\p_i)=0$, then the functions $S(\p_i)$
satisfy the same relation
$$
\sum_i r_i S(\p_i)=0.
$$
\end{proposition}
\begin{example}\label{Scone}.

\noindent$\bullet$ If $\p=\{s\}$ is a point,  then we have two
cases. If $s$ is an integral point, then $S(\{s\})(\xi)=e^{\la
\xi,s\ra}$, otherwise $S(\{s\})=0$.

\noindent$\bullet$  In dimension $1$, if  $\p=s+\R^+$, where $s\in
\Q$, then
$$
S(s+\R_+)(\xi) =\frac{e^{k\xi}}{1-e^{\xi}}
$$
where  $k$ is  the smallest integer greater or equal than $s$.

\noindent$\bullet$ Let $\a$ be a solid simplicial affine cone with
vertex $s\in V_\Q$. Let $ v_1, \dots, v_d$ be integral generators of
the edges of $\a$. Then
\begin{equation}\label{Sparallelepiped} S(\a)(\xi)=
\left(\sum_{x\in (s+\Box( v_1, \dots, v_d))\cap \lattice} e^{\la
  \xi,x\ra} \right)\prod_{i=1}^d \frac{1}{1- e^{\la \xi,v_i\ra}}.
\end{equation}
We obtain this formula by observing that any element $x$ of the
affine cone $\a$ can be written in a unique way as a  sum $ y +
\sum_{i=1}^d n_i v_i $ where $y$ lies in the semi-open
parallelepiped  $s+ \Box( v_1, \dots, v_d)$ and the coefficients
$n_i$ are non negative  integers, and that the point $x$ is integral
if and only if $y$ is.

Thus indeed  $S(s+\c)$ can be extended to a meromorphic function on
the whole of $V^*$.
\end{example}

Let us check the valuation property in dimension one:
\begin{example}\label{one}
In dimension one, let $s\in \Q$.
The relation $\chi(s+\R_+)+\chi(s-\R_+)-\chi(\{s\})=\chi(\R)$ must imply
 $S(s+\R_+)+S(s-\R_+)=S(\{s\})$. If $s=0$, we have indeed
$$
\frac{1}{1-e^{\xi}}+\frac{1}{1-e^{-\xi}} = 1=S(\{s\})(\xi),
$$
while, if  $0<s<1$, then
$$S(s + \R_+)(\xi)= \frac{e^\xi}{1-e^\xi},$$
$$S(s-\R_+)(\xi)=\frac{1}{1-e^{-\xi}},$$
thus $S(s+\R_+)+S(s-\R_+)= 0 = S(\{s \}).$
\end{example}
The valuation property of the maps  $\p\to S(\p)$ and $\p\to
I(\p)$ have the following important corollary. This  was first
obtained by Brion \cite{b} using toric varieties:
\begin{theorem}(Brion)
Let $\p$ be a polyhedron in $V$. Then
\begin{eqnarray*}
I(\p) &= &\sum_{v\in\CV(\p)}I(\t(\p,v)),\\
S(\p) &=& \sum_{v\in\CV(\p)}S(\t(\p,v)).
\end{eqnarray*}
\end{theorem}
The singularities of the functions $I(\a)$ and $S(\a)$ are easy to
compute for a pointed affine  cone:
\begin{lemma}\label{rdelta}
Let $\a= s+\c$ be a pointed affine  cone with vertex $s$ and let
$v_1, \dots, v_k$  be rational generators of the  edges of the
cone $\c$. The products
$$
\left(\prod_{i=1}^k \la \xi,v_i\ra \right)I(s+\c)(\xi) \mbox{ and}
\left(\prod_{i=1}^k \la \xi,v_i\ra \right)S(s+\c)(\xi)
$$
are  analytic near $0$.
\end{lemma}
\begin{proof}
It is easy to see that the cone $\c$ admits a subdivision into
simplicial cones whose edges are already edges of $\c$. Thus,
thanks to  the valuation properties of Proposition
\ref{valuationS}, it is enough to prove the lemma when $\c$ itself
is a simplicial cone. In this case it follows immediately from
Formulas (\ref{Iparallelepiped}) and (\ref{Sparallelepiped}).
$\Box$
\end{proof}

\section{The main construction}\label{mainconstruction}

In this section, we will perform the main construction  of this
article: to any affine cone $\a $ in $V$ or in a rational quotient
$V/L$ of $V$, we will associate an \textbf{analytic} function
$\mu(\a)$  defined in a neighborhood of $0$ in $V^*$. In the next
section, if $\p$ is a convex rational polytope in $V$, we will
obtain a local Euler-Maclaurin formula for $\p$ in terms  of the
functions $\mu(\t(\p,\f))$ associated to the transverse cones
$\t(\p,\f)$  of $\p$ along its various faces $\f$.

 We will denote the ring  of analytic functions with rational coefficients, defined
in a neighborhood of $0$ in $V^*$,    by $\CH(V^*)$ and  the ring of
meromorphic functions with rational coefficients, defined in a
neighborhood of $0$ in $V^*$,    by $\CM(V^*)$

 We  will  need to
extend to the  space $V^*$ some meromorphic functions which are a
priori defined only on a subspace of the form $(V/L)^*= L^\perp $ of
$V^*$.  For that purpose, we fix a scalar product $Q(x,y)$ on $V$.
We assume that $Q$ is rational, meaning that  $Q(x,y)$ is rational
for $x,y \in V_{\Q}$. We denote also by $Q(\xi,\eta)$ the dual
scalar product on $V^*$  and we use the orthogonal projection
$\proj_{ L^\perp}: V^*\to L^\perp$. If $\phi$ is a meromorphic
function (with rational coefficients) on $L^\perp$, we still write
$\phi$ for the function on $V^*$ defined by $\xi\mapsto
\phi(\proj_{L^\perp}(\xi))$. It is meromorphic with rational
coefficients.

Actually, we will do this not only for the space $V$ itself, but
also when $V$ is replaced by a rational quotient space $W$. The
dual $W^*$ is a subspace of $V^*$, thus it inherits the scalar
product of $V^*$.

\noindent {\bf A word of caution.} {\em  Let $L$ be a rational
subspace of $V$. By means of the scalar product $Q$, we can
identify $V/L$ with $L^{\perp_Q}\subset V$, the orthogonal of $L$
with respect to $Q$. However these two spaces  are not isomorphic
as rational spaces. The lattice of $V/L$ corresponds to the
orthogonal projection of $\lattice$ on   $L^{\perp_Q}$;   it
contains the lattice $L^{\perp_Q}\cap \lattice$, and the inclusion
is strict in general, see Figure \ref{normal}.}

\medskip
\begin{figure}[!ht]
\begin{center}
\includegraphics{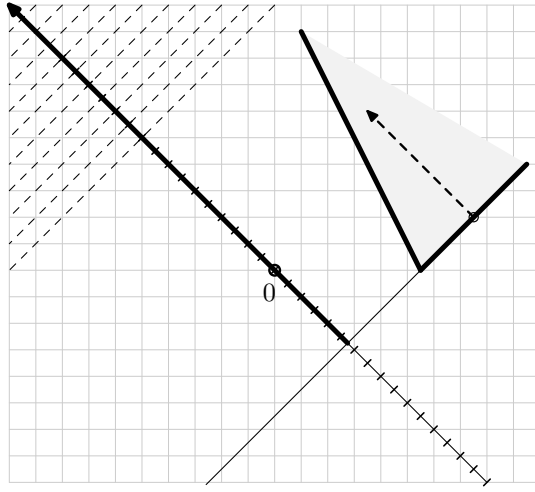}
\caption{In dimension $2$, the transverse cone along an edge with
its lattice}\label{normal}
\end{center}
\end{figure}
Let $\a$ be an affine  cone in $V$ and let $\f$ be a face of $\a$.
Recall that  the transverse cone $\t(\a,\f)$ is the projection
$\pi_\f(\a)$ of $\a$ in $V/\lin(\f)$. When we identify $V/\lin(\f)$
with the orthogonal $\lin(\f)^{\perp_Q}$, the transverse cone
$\t(\a,\f)$ is a  pointed affine cone in $\lin(\f)^{\perp_Q}$.

\bigskip

We denote by ${\CC_{pointed}(V)}$ the set of pointed affine cones in
$V$.

\begin{proposition}\label{defdemu}
Let $V$ be a rational space and let $Q$ be a rational scalar product
on $V$. There exists a unique family of maps  $\mu_W$, indexed by
the rational quotient spaces $W$   of $V$,  such that the family
enjoys the following properties.

(a)  $\mu_W$ maps  $\CC_{pointed}(W)$ to  $\CM(W^*)$.

(b) If $W=\{0\}$,  then  $\mu_W(\{0\})= 1 $.

(c) For any pointed affine cone $\a$ in $W$, one has
\begin{equation}\label{definitionmu}
S(\a)= \sum_{\f\in\CF(\a)}\mu_{W/\lin(\f)}(\t(\a,\f))I(\f)
\end{equation}
where the sum is over the set of  faces of $\a$.
\end{proposition}
In Formula (\ref{definitionmu}), the transverse cone $\t(\a,\f)$ is
a pointed affine cone in the quotient space $W/\lin(\f)$. The
function $\mu_{W/\lin(\f)}(\t(\a,\f))$ is a meromorphic function on
a neighborhood of $0$ in the dual $(W/\lin(\f))^*$. We give a
meaning to the formula by extending  this function to a neighborhood
of $0$ in the whole space $W^*$ by means of orthogonal projection.
The function $I(\f)$ is defined as a meromorphic function on $W^*$
as in Section 3.

\begin{proof}
The result is easily obtained by induction on the dimension of
$W$. If $W=\{0\}$, the only cone is  $\a=\{0\}$ and Formula
(\ref{definitionmu}) is true. Let $\a$ be a pointed affine cone in
$W$. Let $s$ be the vertex of $\a$. The transverse cone at the
zero-dimensional face $s$ is $\a$ itself.  Formula
(\ref{definitionmu}) gives
\begin{equation}\label{defmu-induction}
S(\a)(\xi) = e^{\la \xi,s\ra }\mu_W(\a)(\xi)+ \sum_{\f,
\dim\f>0}\mu_{W/\lin(\f)}(\t(\a,\f))(\xi)I(\f)(\xi).
\end{equation}

For a face $\f$ of positive dimension, the transverse cone
$\t(\a,\f)$ is a  pointed affine cone in the vector space
$W/\lin(\f)$. Therefore,  $\mu_{W/\lin(\f)}(\t(\a,\f))$ being
defined by the induction hypothesis, Formula (\ref{defmu-induction})
defines $\mu_W(\a)$ in a unique way, as a meromorphic function on
$W^*$.$\Box$
\end{proof}

The following property follows immediately from the definition:
\begin{proposition}
If $V_1\subset V_2$ and $\a$ is an affine  cone contained in $V_1$,
then the function $\mu_{V_2}(\a)$  is the lift to $V_2^*$ of
$\mu_{V_1}(\a)$, by the natural restriction map $V_2^*\to V_1^*$ .
\end{proposition}
In the rest of this article, we will omit the subscript $W$ in the
notation $\mu_W(\a)$.
\begin{proposition}
The functions  defined in Proposition \ref{defdemu} have the
following properties:

(a) For any $s\in \lattice$, one has $\mu(s +\a)= \mu(\a)$.

(b) The map $\a\mapsto \mu(\a)$  is equivariant with respect to
lattice-preserving  linear isometries. In other words, let $g$ be
a linear isometry of $W$ which preserves the lattice $\lattice$
and denote its transpose  by $^{t}g$, then $\mu(g(\a))(
^{t}g^{-1}\xi)=\mu(\a)(\xi)$.

(c) The map $\a\mapsto \mu(\a)$ is multiplicative with respect to
orthogonal sums of cones. More precisely, if $W$ is an orthogonal
sum $W=W_1\oplus W_2$ and  $\a_i$ is an affine cone in $W_i$ for
$i=1,2$, then
$$
\mu(\a_1 + \a_2)= \mu(\a_1)\mu(\a_2).
$$

(d)  If  $\a\in \CC_{pointed}(W)$ is such that   $<\a>\cap
\lattice_W =\emptyset$, then
 $\mu_W(\a)=0$.
\end{proposition}
\begin{proof} The invariance in {\it (a)} and {\it (b)}  and the multiplication rule  in {\it (c)}
follow immediately from the definition, by induction. To prove {\it
(d)} , assume that $<\a>\cap \lattice_V =\emptyset$. Then $S(\a)=0$
and, for any face $\f$ of $\a$, the transverse cone $\t(\a,\f)$ does
not contain any integral point of  $V/\lin(\f)$, therefore, by
induction, $\mu(\t(\a,\f))=0$ for $\dim \f>0$, hence $\mu(\a)=0$.
$\Box$
 \end{proof} \medskip

If $\a$ is an affine cone in $V$ which contains a straight line,
we define
 $\mu(\a)=0$. Since all faces $\f$ of $\a$ contain a straight line, Formula
 (\ref{definitionmu}) still holds in this case. Thus we have defined
$\mu(\a)$ for any rational affine cone in any quotient space of
$V$.

\bigskip

Our objective  is to show that $\mu(\a)$ is indeed analytic near $0$, but we
will first observe some further properties of this family of functions.

It is easy to  compute $\mu$ in dimension $1$. Let $t\in \Q$ such that
 $0 \leq t <1$. Then $\mu(- t+\R_+)(\xi)$ is defined by
$$
S(-t+\R_+)(\xi) = \frac{1}{1-e^\xi}=  e^{-t\xi}\mu(-t+\R_+)(\xi) +
\int_{-t}^\infty e^{x\xi}dx
$$
hence
\begin{equation}\label{mudim1}
\mu(- t+\R_+)(\xi)=\frac{e^{t\xi}}{1-e^\xi}+\frac{1}{\xi}.
\end{equation}
We may write this in terms of the Bernoulli polynomials $b(n,t)$,
defined by the generating series (\ref{bernoulli}). We obtain
\begin{equation}\label{mudim1bern}
\mu(- t+\R_+)(\xi)= - \sum_{n=0}^\infty \frac{b((n+1),t)}{(n+1)!}\xi^n.
\end{equation}

Let $\a$ be a 1-dimensional pointed cone in $V $, that is to say a
half-line. If $\a$ does not contain any integral point, then
$\mu(\a)=0$. If $\a$ contains integral points, then there exists an
integral point $a \in \a$ such that the translated half-line $ \a-a$
is of the form
$$
(-t+\R_+)v
$$
 where  $v$ is a primitive integral vector and $t \in [0,1[$.

By a similar computation, we get, for $\xi \in V^*$,
\begin{equation}\label{mudim1bis}
\mu(\a)(\xi)= \frac{e^{t
\la\xi,v\ra}}{1-e^{\la\xi,v\ra}}+\frac{1}{\la\xi,v\ra}.
\end{equation}

\medskip
The next step is crucial to our construction; we  will prove that
the map $\c\to \mu(s+\c)$ enjoys the valuation property:
\begin{proposition}
Let $\c_i$ be a finite family of cones in $V$. Assume that there
exists a linear relation between their characteristic functions
$\sum _i r_i \chi(\c_i)=0.$ Then, for any $s\in V_\Q$, we have the
corresponding relation $\sum _i r_i \mu(s+\c_i)=0.$
\end{proposition}
\begin{example} Consider the subdivision of Example \ref{subdivisionR}
  in dimension one.  Let  $s\in\Q$. Then
$$
0= \mu(\R)=\mu(s+\R_+)+\mu(s-\R_+)-\mu(\{s\}).
$$
Indeed, if $s$ is an integer, we have
$$
\mu(s+\R_+)(\xi)+\mu(s-\R_+)(\xi)= (\frac{1}{1-e^\xi}+
\frac{1}{\xi})\; + (\frac{1}{1-e^{-\xi}}-\frac{1}{\xi})= 1=
\mu(\{s\})(\xi)
$$
while, if $s \in   ]-1,0[$, we have
  $$\mu(s+\R_+)(\xi)+\mu(s-\R_+)(\xi)=(\frac{e^{-s\xi}}{1-e^\xi}+\frac{1}{\xi}) \; +
(\frac{e^{-(s+1)\xi}}{1-e^{-\xi}}-\frac{1}{\xi})=0=
\mu(\{s\})(\xi).$$
\end{example}
\begin{proof}

We will prove the proposition  by induction on $\dim V$. By a
standard  argument (\cite{Sallee}, see also \cite{brionbki}), it
suffices to prove the result in  the following particular case
(see Figure \ref{cut}).
\begin{figure}[!h]
\begin{center}
\includegraphics{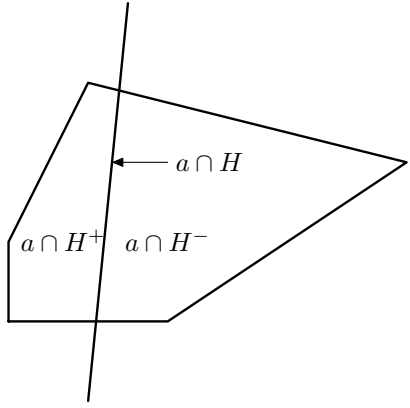}
\caption{The cut of a 3-dimensional cone with 5 edges by a 2-dimensional plane H.
The cone is represented by its slice in the figure plane.
 Thus $H$ is represented by a line.}\label{cut}
\end{center}
\end{figure}
Let $\a$ be a pointed solid affine cone in $V$ with vertex $s$,
let $H$ be  an affine  hyperplane through $s$. Denote by $H^\pm$
the closed half-spaces separated by $H$. Then we have
$$
\chi(\a)= \chi(\a \cap H^+)+\chi(\a\cap H^-)-\chi(\a\cap H)
$$
and we must prove:
\begin{equation}\label{coupure-hyperplan}
\mu(\a)= \mu(\a \cap H^+)+\mu(\a\cap H^-)-\mu(\a\cap H).
\end{equation}

We proceed to prove (\ref{coupure-hyperplan}).

The functions $S(\a)$ have the valuation property
$$
S(\a)-S(\a \cap H^+)-S(\a\cap H^-)+S(\a\cap H)=0.
$$
By applying Formula (\ref{definitionmu}), we obtain
the following expansion of the left hand side
\begin{eqnarray}\label{expansion}
&& \sum_{\f\in\CF(\a)}\mu(\t(\a,\f))I(\f)
 - \sum_{\f\in\CF(\a\cap   H^+)} \mu(\t(\a\cap H^+,\f))I(\f)\\
\nonumber
&& - \sum_{\f\in\CF(\a\cap
  H^-)}\mu(\t(\a \cap H^-,\f))I(\f)
+ \sum_{\f\in\CF(\a\cap H)}\mu(\t(\a\cap H,\f))I(\f).
\end{eqnarray}
Let $L$ be a affine subspace of $V$ of dimension $>0$. We will show
\begin{eqnarray}\label{expansionL}
&& \sum_{\f\in\CF(\a),<\f>= L}\mu(\t(\a,\f))I(\f)
 - \sum_{\f\in\CF(\a\cap   H^+), <\f>= L} \mu(\t(\a\cap H^+,\f))I(\f)\\
\nonumber
&& - \sum_{\f\in\CF(\a\cap
  H^-), <\f>= L}\mu(\t(\a \cap H^-,\f))I(\f)
+ \sum_{\f\in\CF(\a\cap H), <\f>= L}\mu(\t(\a\cap H,\f))I(\f)\\
\nonumber &=&0.
\end{eqnarray}

From the relation (\ref{expansionL}), it follows that the terms in
(\ref{expansion}) corresponding to the faces $\f$ of positive
dimension add up to $0$. Therefore  the contribution of the
$0$-dimensional vertex $\{s\}$  to (\ref{expansion}) is also equal
to $0$, which  proves the relation (\ref{coupure-hyperplan}).

We fix $L$ and we proceed to prove   (\ref{expansionL}).

Remark that all the transverse cones which appear in
(\ref{expansionL}) are affine cones in $V/\lin(L)$. We will apply
the induction hypothesis to $V/\lin(L)$.

 I) First we consider the case where there is a face $\f$ of $\a$ such
that $<\f>= L$. There are three cases, according to whether the
relative interior of $\f$ meets both the interiors of  $H^\pm$, or
only one, or none of them (in the third case,  $\f$ is contained
in $H$).
\begin{figure}[!h]
\begin{center}
\includegraphics{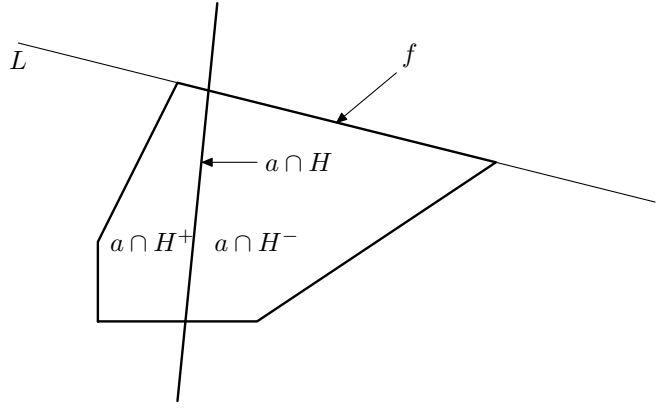}
\caption{Case I.1}\label{cut1}
\end{center}
\end{figure}

$\bullet$ Case I.1: The relative interior of
$\f$ meets both the interiors of  $H^\pm$. Then $\f\cap H^\pm$ is a face
of  $\a \cap H^\pm$ and  $<\f\cap H^\pm>= L$.
Thus we have to prove
$$
\mu(\t(\a,\f))I(\f)- \mu(\t(\a \cap H^+,\f\cap H^+))I(\f\cap H^+ ) -
\mu(\t(\a \cap H^-,\f\cap H^-)) I(\f\cap H^-)=0.
$$
The three transverse cones $\t(\a,\f)$, $\t(\a \cap H^+,\f\cap H^+)$
and  $\t(\a \cap H^-,\f\cap H^-)$ coincide. The integrals add up:
$$
I(\f)= I(\f\cap H^+ )+I(\f\cap H^-),
$$
thus we get
$$
\mu(\t(\a,\f))\left(I(\f)-I(\f\cap H^+ )-I(\f\cap H^-)\right)
$$
which is equal to $0$ as required.

\begin{figure}[!ht]
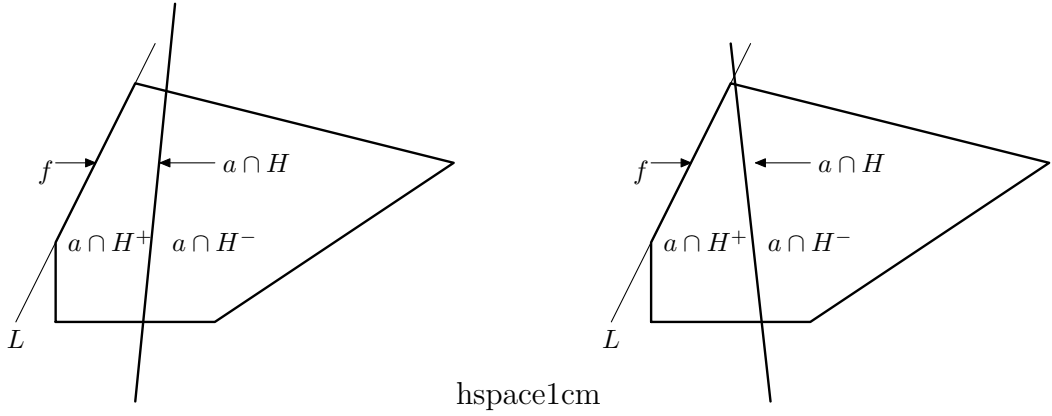

\begin{center}
\includegraphics{proofsNicole.04.eps}hspace{1cm}\includegraphics{proofsNicole.41.eps}
\caption{Case I.2}\label{cut2}
\end{center}
\end{figure}
$\bullet$ Case I.2: The relative interior of $\f$ is contained in
the interior of, say,  $H^+$.   Then $\f= L\cap\a\cap H^+$ is also
a face  of $\a\cap H^+$,   but  $L\cap\a \cap H^-$ and
$L\cap\a\cap H$ are smaller dimensional, or empty. This time we
have to prove
$$
\mu(\t(\a,\f))I(\f)- \mu(\t(\a \cap H^+,\f\cap H^+))I(\f\cap H^+
)=0.
$$
The transverse cones $\t(\a,\f)$ and $\t(\a \cap H^+,\f)$ coincide,
so  we get
$$
\mu(\t(\a,\f))\left(I(\f)-I(\f\cap H^+ )\right).
$$
As $\f= \f\cap H^+ $, this is equal to $0$ as required.

\begin{figure}[!h]
\begin{center}
\includegraphics{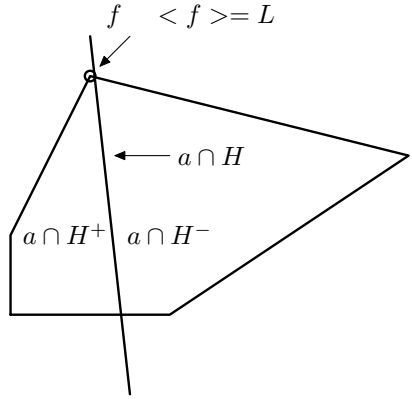}
\end{center}
\caption{Case I.3}\label{cut3}
\end{figure}

$\bullet$ Case I.3: $\f$ is contained in $H$, thus it is a face of all four
cones, in other words
$$
\f=\f\cap H^+= \f\cap H^-= \f\cap H.
$$
This time we have to prove
\begin{eqnarray*}
& &\mu (\t(\a,\f))I(\f) - \mu(\t(\a \cap H^+,\f\cap H^+))I(\f\cap H^+ )\\
&- & \mu (\t(\a \cap H^-,\f\cap H^-)) I(\f\cap H^-) + \mu(\t(\a\cap
H,\f))I(\f\cap H)=0.
\end{eqnarray*}
In this case, the intersection of the transverse cones is
$$
 \t(\a \cap H^+,\f)\cap \t(\a \cap H^-,\f)= \t(\a\cap H,\f),
$$
the union of the transverse cones is
$$
\t(\a,\f)= \t(\a \cap H^+,\f)\cup \t(\a \cap H^-,\f).
$$
Thus we get
$$
\Big(\mu(\t(\a,\f)) - \mu(\t(\a \cap H^+,\f)) - \mu( \t(\a \cap
H^-,\f))+ \mu(\t(\a\cap H,\f))\Big)I(\f).
$$
From the induction hypothesis applied to the space $V/\lin(L)$, we
deduce  that this is equal to $0$.

\begin{figure}[!h]
\begin{center}
\includegraphics{proofsNicole.06.eps}hspace{1cm}\includegraphics{proofsNicole.07.eps}
\caption{Case II }\label{caseII}
\end{center}
\end{figure}
II) Next, we consider an affine  subspace $L$ such that $L\cap \a$ is not a face of
$\a$ but, say,   $\f= L\cap\a\cap H^+$  is a
face of $\a\cap H^+$.  Then we must have $L\subset H$ so that
 $\f$ is a face of the three cones  $\a\cap H$, $\a \cap
 H^+$ and  $\a \cap H^-$, but not a face of $\a$.
We have to show that
\begin{eqnarray*}
&&\mu(\t(\a \cap H^+,\f\cap H^+))I(\f\cap H^+ )
+ \mu (\t(\a \cap H^-,\f\cap H^-)) I(\f\cap H^-) \\
&&-  \mu(\t(\a\cap H,\f\cap H))I(\f\cap H)=0,
\end{eqnarray*}
with
$$
\f\cap H^+=\f\cap H^-=\f\cap H=\f.
$$
In this case, the union $\t(\a \cap H^+,\f)\cup \t(\a \cap H^-,\f)$
is the
 projection $\pi_{\lin(\f)}(\a)$ of $\a$ on $V/\lin(\f)$; it is not
 pointed, therefore, applying again  the induction hypothesis
to the space $V/\lin(L)$,  we have
$$
\mu( \t(\a \cap H^+,\f))+\mu(\t(\a \cap H^-,\f))-\mu( \t(\a\cap
H,\f))=\mu( \pi_{\lin(\f)}(\a))= 0
$$
and the result follows.
$\Box$
 \end{proof}
 \bigskip

\begin{corollary}\label{sommeconesausommet}
Let $\p$ be a polytope in $V$ and $s\in V_\Q$.  Then
$$
\sum_{v\in\CV(\p)}\mu(s+ \dir(\t(\p,v)))
$$
is equal to $1$ if the point $s$ is integral and $0$ otherwise.
\end{corollary}
\begin{proof}
This follows immediately from the valuation property and the
relation (see \cite{barpom}) between the characteristic functions:
$$
\sum_{v\in\CV(\p)}\chi(\dir(\t(\p,v)))=\chi(\{0\}) \mbox{ mod}
\,\,\CL
$$
where $\CL$ denotes the space  of linear combinations of
characteristic functions of cones with lines. $\Box$
\end{proof}

\medskip
Now we show that our functions are analytic near $0$.
\begin{proposition}\label{analytic}
Let $\a$ be an affine cone in $V$. The function  $\mu(\a)$ is
analytic near $0$.
\end{proposition}
\begin{proof}
The result is true when $V=0$ (and the explicit computation shows
that it is true also when $\dim V=1$). We will prove it by
induction on the dimension of $V$. Using  the valuation property,
it is enough to prove the analyticity when $\a$ is a solid
simplicial unimodular affine cone in $V$. Let $v_1,\dots,v_k$ be
primitive integral generators of the edges of $\dir\a$. If $\Phi$
is a meromorphic function on $V^*$ such that the product
$$
\left(\prod_{i=1}^k \la \xi,v_i\ra \right)\Phi(\xi)
$$
is analytic, we denote by  $\Res_{v_1}(\Phi)$ the residue of
$\Phi$ along $v_1=0$, that is to say the  restriction to
$v_1^\perp \subset V^*$ of
$$
\la\xi,v_1\ra \Phi(\xi).
$$
From  the properties of the functions $S$ and $I$ (Lemma
\ref{rdelta}) and the induction hypothesis, it follows that the
product
$$
\left(\prod_{i=1}^k \la \xi,v_i\ra \right)\mu(\a)(\xi)
$$
is analytic. Thus we want to show that $\Res_{v_1}(\mu(\a))=0$.
Starting from  the  defining formula (\ref{definitionmu}), we want
to prove that, for $\xi\in v_1^\perp $, we have
\begin{equation}\label{residu}
\Res_{v_1}(S(\a))(\xi)= \sum_{\f;
\dim\f>0}\mu(\t(\a,\f))(\xi)\Res_{v_1}(I(\f))(\xi).
\end{equation}
Let us denote by $\pi$ the projection $V\to V/<v_1>$.  The cone
$\pi(\a)$ is also a simplicial unimodular cone with primitive
integral generators $\pi(v_2),\dots,\pi(v_d)$. If $\a$ is
unimodular, then the parallelepiped $s+\Box(v_1,\dots,v_d)\subset
V$ contains exactly one integral point. Therefore,  the explicit
computations (\ref{Iparallelepiped}) and (\ref{Sparallelepiped})
of  $I(\a)$ and $ S(\a)$  imply  immediately that
 the residues along $v_1=0$ of the functions
$S(\a)$ and $I(\a)$ are given by:
$$
\Res_{v_1}(S(\a)) = - S(\pi(\a)),
$$
$$
\Res_{v_1}(I(\a)) = - I(\pi(\a)).
$$
In the sum (\ref{residu}), only the faces $\f$ for  which $v_1$ is
an edge of $\dir \f$ contribute, and these faces are in one to one
correspondence with the faces  of $\pi(\a)$. For such a face $\f$,
the transverse cone of $\pi(\a)$ along $\pi(\f)$ coincides with the
transverse cone $\t(\a,\f)$ and we have also
$$
\Res_{v_1}(I(\f)) = - I(\pi(\f)),
$$
whence (\ref{residu}), and the proposition.
$\Box$
\end{proof}
\medskip

Next we will show that Formula (\ref{definitionmu}) still holds
when $\a$ is replaced by any polyhedron $\p$. This will be an easy
consequence of Brion's theorem and the valuation property of $I$
and $S$.  In the following three theorems, we  collect the results
of this section.

\begin{theorem}\label{maintheorem}

Let $V$ be a rational space and $Q$ a
rational scalar product on $V^*$. If $W=V/L$ is a rational quotient space
 of $V$,
we denote by $\CC(W)$ the set of affine cones in $W$. For $\a\in
\CC(W)$, let $I(\a)$ and $S(\a)$ be the meromorphic functions with
rational coefficients on $W^*$ defined in Propositions
\ref{valuationI} and \ref{valuationS}.

There exists a unique family of maps  $\mu_W$, indexed by the
rational quotient spaces $W$   of $V$, such that the family enjoys
the following properties:

(a)  $\mu_W$ maps  $\CC(W)$ to  $\CH(W^*)$, the space of analytic
functions on $W^*$,  with rational Taylor coefficients.

(b) If $W=\{0\}$,  then  $\mu_W(\{0\})= 1 $.

(c) If the affine cone $\a\in \CC(W)$ contains a straight line, then
$\mu_W(\a)=0.$

(d) For any affine cone $\a$ in $W$, one has
$$
S(\a)= \sum_{\f\in\CF(\a)}\mu_{W/\lin(\f)}(\t(\a,\f))I(\f)
$$
where the sum is over all faces of the cone $\a$.

\end{theorem}
As in all this section, in  Formula {\it (d)}, the function
$\mu_{W/\lin(\f)}(\t(\a,\f))$ is considered as a function on $W^* $
itself by means of the orthogonal projection $W^* \to
(W/\lin(\f))^*= (\lin(\f))^\perp $ with respect to the scalar
product on $W^*\subset V^*$.

\begin{theorem}\label{maintheoremplus}
The analytic functions  defined in Theorem \ref{maintheorem} have
the following properties:

(a) For any $x\in \lattice$, one has $\mu(x +\a)= \mu(\a)$.

(b) The map $\a\mapsto \mu(\a)$  is equivariant with respect to
lattice-preserving  isometries. In other words, let $g$ be an
isometry of $W$ which preserves the lattice $\lattice$. Then
$\mu(g(\a))( ^{t}g^{-1}\xi)=\mu(\a)(\xi)$.

(c) If $W$ is an orthogonal sum $W=W_1\oplus W_2$ and  $\a_i$ is
an affine cone in $W_i$ for $i=1,2$, then
$$
\mu(\a_1 + \a_2)= \mu(\a_1)\mu(\a_2).
$$

(d)  For a fixed $s\in W_\Q$, the map $\c\to \mu(s+\c)$ is a
valuation on the set of cones in $W$.

(e)  Let $\p\subset W$ be a polyhedron, then
\begin{equation}
S(\p)(\xi)= \sum_{\f\in\CF(\p)}\mu(\t(\p,\f))(\xi)I(\f)(\xi).
\end{equation}
\end{theorem}
\begin{proof}
In Theorem \ref{maintheoremplus}, only point {\it (e)}  has not yet
been proven. If $v$ is a vertex of $\p$, let us denote by
$\CF(\p,v)$ the set of faces of $\p$ which contain $v$. For such a
face $\f\in\CF(\p,v) $, the intersection $\d=<\f>\cap \t(\p,v)$ is a
face of the cone $ \t(\p,v)$ and this correspondence is a bijection
between $\CF(\p,v)$ and $\CF(\t(\p,v))$ with inverse given by $\f =
\d\cap \p$. Moreover the transverse cone  $ \t(\p,v)$ of $\p$ along
its face $<\f>\cap \t(\p,v)$ coincides with the transverse cone $
\t(\p,\f)$ of $\p$ along $\f$. Therefore we have
$$
 S(\n (\p,v))= \sum_{\f\in \CF(\p,v)}\mu(\t(\p,\f))\; I(<\f>\cap \t(\p,v)).
$$
Replacing $S(\t(\p,v))$ with the right-hand side of this equality in
Brion's formula, we obtain
$$
S(\p)= \sum_{v\in \CV(\p)}\sum_{\f\in \CF(\p,v)}\mu(\t(\p,\f))\;
I(<\f>\cap \t(\p,v)).
$$
Then we reverse the order of summation and get
$$
S(\p)= \sum_{\f\in
\CF(\p)}\mu(\t(\p,\f))\sum_{v\in\CV(\f)}I(<\f>\cap \t(\p,v)).
$$
The last sum is equal to $I(\f)$. $\Box$  \end{proof}

\begin{theorem}\label{calculable}
Assume $V=\R^d$ with a fixed dimension $d$. Then, for $m$ fixed,
there exists a polynomial time algorithm which computes $\mu(\a)$
at order $m$  for any affine cone $\a\subset V$.
\end{theorem}
\begin{proof}
By  \cite{barvinok}, there exist polynomial time algorithms which
compute  the functions $I(\c)$ and $S(\c)$ at order $m$ for any
affine cone $\c$ in $\R^k$, if $k\leq d$. Therefore by induction
we get an algorithm which computes $\mu(\a)$ for any $\a\subset
V$. $\Box$  \end{proof}

\medskip Let $\sigma$ be a cone in the dual space $V^*$. The dual
cone $\sigma^*\subset V $ contains the vector subspace
$<\sigma>^\perp$. Let us denote by $\pi_{<\sigma>^\perp}$ the
projection $V \to V/<\sigma>^\perp$.  For any $s\in V_\Q$, the
projected cone  $\pi_{<\sigma>^\perp}(s+\sigma^*)$
 is a pointed cone in $V/<\sigma>^\perp$. Thus
 $\mu(\pi_{<\sigma>^\perp}(s+\sigma^*))$ is  an analytic function  on
$(V/<\sigma>^\perp)^* \cong  <\sigma>\subset V^*$. We consider it as a
 function
 on $V^*$ by means of orthogonal projection, as before. We
 obtain a map  $\mu_s^*: \CC(V^*) \to \CH(V^*)$ defined by:
\begin{definition}\label{mustar}
$$
\mu_s^*(\sigma)= \mu(\pi_{<\sigma>^\perp}(s+\sigma^*)).
$$
\end{definition}
From the  valuation behavior  of  $\mu$, it follows that $\mu_s^*$
is a solid valuation. In other words, the following corollary holds.
\begin{corollary}
Let $\sigma$ be a cone in $V^*$, and let $\{\sigma'\}$ be a
subdivision of $\sigma$. For any  $s\in V_\Q$, we have
$$
\mu_s^*(\sigma)= \sum_{\dim \sigma'= \dim \sigma}
\mu_s^*(\sigma').
$$
\end{corollary}
\begin{proof}
 As  $\{\sigma'\}$ is  a subdivision of $\sigma$, we have
$$
\chi(\sigma)= \sum_{\dim \sigma'= \dim \sigma}\chi(\sigma') +
\sum_{\dim \sigma'< \dim \sigma} \pm \chi(\sigma').
$$
Let $L=<\sigma>^\perp\subset V$ and let $\pi$ denote the
projection $V\to V/L$. The map $\tau \mapsto \chi(\pi(s+\tau^*))$
is a valuation on the set of cones in $V^*$ (see \cite{barpom} for
instance). Therefore
$$
\chi(\pi(s+\sigma^*))= \sum_{\dim \sigma'= \dim \sigma}\chi(\pi(s+\sigma'^*)) +
\sum_{\dim \sigma'< \dim \sigma} \pm \chi(\pi(s+\sigma'^*)).
$$
If $\dim \sigma'<\dim \sigma$,  then the cone
$\pi(s+\sigma'^*)\subset V/L$ contains a straight line , thus
$\mu(\pi(s+\sigma'^*))=0$, and the corollary follows from the
valuation property of $\mu$. $\Box$\end{proof}

In a companion paper \cite{be-ve-todd}, we will prove the following
theorem, which extends to equivariant homology a result of
Pommersheim-Thomas \cite{pomm} by which they answered a question of
Danilov \cite{danilov}.
\begin{theorem}\label{todd}
Let $\CE$ be a fan in $V^*$ and let $X$ be the corresponding toric
variety. For $\sigma\in \CE$, let $X(\sigma)\subseteq X$ be the
corresponding orbit closure. Then the equivariant Todd class of $X$
is equal to
$$
\sum_{\sigma\in \CE}\mu_0^*(\sigma)[X(\sigma)]
$$
in the equivariant homology ring of $X$.
\end{theorem}
\section{Local Euler-Maclaurin formula}

As in the previous section,  $V$ is a rational space and we fix a
scalar product on $V$.  Let  $\p$ be a (convex rational)
polyhedron in $V$. To each face $\f$ of $\p$, we are going to
associate a linear differential operator $D(\p,\f)$ on $V$.

To any analytic function $\Phi(\xi)$ on  $V^*$, defined near $0$,
there corresponds a unique linear differential operator $D(\Phi)$
(of infinite degree) with constant coefficients on $V$ such that
$\Phi(\xi)$ is the symbol of $D(\Phi)$. More precisely, for
$\xi\in V^*$, let us denote by $e^\xi$ the function  $x\mapsto
e^{\la\xi,x\ra}$ on $V$, then  $D(\Phi)$ is defined by the
relation
$$
D(\Phi)\cdot e^\xi= \Phi(\xi)e^\xi \mbox{  for } \xi \mbox{ small
enough}.
$$
Let $W=V/L $ be a quotient space of $V$ and let  $\a$ be  a
pointed affine cone in $W$.  In the previous section, we
constructed an  analytic function $\mu(\a)$ on $W^*= L^\perp
\subset V^*$. By orthogonal projection, we consider $ \mu(\a)$ as
a  function on $V^*$ and we introduce  the corresponding
differential operator $D(\mu(\a))$ on $V$:
\begin{equation}\label{da}
D(\mu(\a))\cdot e^\xi= \mu(\a)(\xi)e^\xi.
\end{equation}

Let $\p$ be a polyhedron in $V$.
\begin{definition}\label{dpf}
Let $\f$ be a face of $\p$. We denote by
$$
D(\p,\f)=  D(\mu(\t(\p,\f)))
$$
the differential operator on  $V$ associated to the transverse cone
$\t(\p,\f)$ of $\p$ along $\f$. We denote its constant term by
$\nu(\p,\f)$. Thus
$$
\nu(\p,\f)= D(\p,\f)\cdot 1=  \mu(\t(\p,\f))(0).
$$
\end{definition}

The operator $D(\p,\f)$, as well as its constant term $\nu(\p,\f)$,
are  {\bf local} in the sense that they depend only on the class of
$ \t(\p,\f)$ modulo integral translations. In particular, if $\p$
has integral vertices, then  $D(\p,\f)$ depends only on the cone of
transverse feasible directions at a generic point of $\f$. The
operator $ D(\p,\f)$ involves only derivatives in directions
orthogonal to the face $\f$.

We are now ready to state the local Euler-Maclaurin formula for any
polytope.
\begin{theorem}(Local Euler-Maclaurin formula)\label{LEMLF}

Let $\p$ be a polytope in $V$. For any polynomial function $h(x)$
on $V$, we have
\begin{equation}\label{maclaurin}
\sum_{x\in \p\cap\lattice}h(x)=\sum_{\f\in\CF(\p)}\int_\f
D(\p,\f)\cdot h
\end{equation}
 where the integral on the face $\f$ is taken with respect to the
Lebesgue measure on $<\f>$ defined by the lattice $\lattice\cap
\lin(\f)$.
\end{theorem}
\begin{proof}

 The method is to check equality (\ref{maclaurin}) for a polynomial of the
 form
$h(x)=\langle \xi,x\rangle ^k$. Taking Taylor series, we may replace
$h(x)$ by $e^{t\langle \xi,x\rangle}$ with $t$ small. Then the
equality (\ref{maclaurin}) becomes the formula in Theorem
\ref{maintheoremplus}, {\it (c)}
$$
S(\p)(\xi)=\sum_{\f\in\CF(\p)} \mu(\t(\p,\f))(\xi)I_\f(\xi).
$$
$\Box$  \end{proof}

In dimension $1$, when $\p$ is an interval $[a_1,a_2]$,  applying
Formulas (\ref{mudim1bis}) and (\ref{mudim1bern}) for
$\mu(a_1+\R^+)$ and $\mu(a_2 - \R^+)$, we obtain
\begin{eqnarray*}
  \sum_{a_1\leq x\leq a_2, x\in \Z}h(x) =
\int_{a_1}^{a_2}h(t)dt
& - & \sum_{n\geq 0}\frac{b((n+1),t_1)}{(n+1)!}((\frac{d}{dt})^n
h)(a_1) \\
& -& \sum_{n\geq 0}(-1)^n
\frac{b((n+1),t_2)}{(n+1)!}((\frac{d}{dt})^n h)(a_2),
\end{eqnarray*}
where $t_1$ and $t_2$ in $[0, 1[ $ are defined by $t_1=k_1 -a_1$,
    $t_2=a_2-k_2$, with $k_1$  the smallest integer greater or equal than
    $a_1$ and $k_2$ the largest integer smaller or equal  than $a_2$
     (Figure \ref{interval}). Of course,
    when $a_1$  and $a_2$ are integers, we recover the historical Euler-Maclaurin formula.

\begin{figure}[!ht]\label{interval}
\begin{center}
\includegraphics{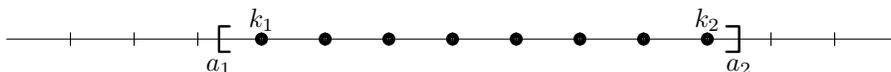}
\caption{Euler-Maclaurin for the interval $[a_1,a_2]$}
\end{center}
\end{figure}
\section{Ehrhart polynomial}
Let $\p$ be a rational polytope in a $d$-dimensional rational
space $V$ and let $q$ be an integer such that $q\p$ has integral
vertices. Let $h(x)$ be a polynomial function of degree $r$ on
$V$. For any integer $t$, we consider the dilated polytope $t\p$
and the corresponding sum
$$
S(t\p,h)= \sum_{x\in t\p\cap\lattice}h(x).
$$
As a function of $t$, it is given by a quasipolynomial: there
exist functions $t\mapsto E_i(\p,h,t)$ on $\Z$ which are  periodic
with period $q$ such that
\begin{equation}\label{coeffEhrhart}
S(t\p,h)= \sum_{i=0}^{d+r} E_i(\p,h,t)t ^i
\end{equation}
whenever $t$ is a positive integer and even in a slightly larger
range including negative values.

\begin{definition}\label{defcoeffEhrhart}
The periodic functions  $E_i(\p,h,t)$ defined by Equation
(\ref{coeffEhrhart}) are called the \emph{Ehrhart coefficients}
\emph{for the polytope} $\p$ \emph{and the polynomial} $h$.

When $h$ is the constant polynomial $h(x)=1$, we denote the
Ehrhart coefficients simply by  $E_i(\p, t)$.
\end{definition}

Our local Euler-Maclaurin formula gives  an expression of the
coefficients $E_i(\p,h,t)$ in terms of the functions
$\mu((\t(\p,\f))$, as we will now explain.

Let $\a$ be a pointed affine cone with vertex $s$  in $V$. If $t$ is
a nonzero integer, we define $\mu(\a, t) = \mu(t\a)$. For $t=0$, we
define $\mu(\a, 0)= \mu(\dir(\a))$. Then $\mu(\a, t+q)= \mu(\a, t) $
for any integer $q\in \N$ such that the point $qs$ is integral. Let
$\p$ be a polyhedron and $\f$ a face of codimension $m$. We define
$D(\p,\f,t)= D(\mu(\t(\p,\f),t))$ for any integer $t\in \N$.
\begin{remark}\label{integralface}
If the affine span $<\f>$ of the face $\f$  contains an integral
point, then the vertex of the transverse cone $\t(\p,\f)$ is
integral, therefore $\mu(\t(\p,\f),t)= \mu(\dir(\t(\p,\f)))$ and
$D(\p,\f,t)$ do not depend on $t$. We have  $D(\p,\f,t)=D(\p,\f)$.
\end{remark}

Let $w_1,\dots,w_m$ be an integral  basis of the subspace
$\lin(\f)^{\perp_{Q}}\subset V$. The operator $D(\p,\f,t)$ has the
following expression:
$$
D(\p,\f,t) = \nu_0(\p,\f,t)+ \sum_{A,|A|=1}^\infty \nu_A
(\p,\f,t)\partial^A
$$
where $A=(a_1,\dots,a_m)$, with $a_i\in\N$ and $\partial^A=
D(w_1)^{a_1}\cdots D(w_m)^{a_m}$. The coefficients $\nu_A
(\p,\f,t)\in \Q$ are periodic with respect to $t$, with period
equal to the smallest integer $q_\f$ such that $q_\f<\f>$ contains
an integral point.
\begin{proposition}
Let $\p$ be a rational polytope and $h$ a polynomial function of
degree $r$ on $V$. Then, for any integer $t\geq 0$, we have
\begin{equation}\label{Stph}
S(t\p,h)=\sum_{\f\in\CF(\p)}\int_{t\f}D(\p,\f,t)\cdot h.
\end{equation}
Furthermore we have
$$
\int_{t\f}D(\p,\f,t)\cdot h = \sum_{i=\dim\f}^{\dim\f+
r}E_i(\p,h,\f,t)\,t^i
$$
where the coefficients $E_i(\p,h,\f,t)$ are periodic with period
$q_\f$.

Hence the Ehrhart coefficients are given by
$$
E_i(\p,h,t)= \sum_{\f,{\rm dim \f} \leq i}E_i(\p,h,\f,t).
$$
\end{proposition}
\begin{proof}
For $t>0$, Formula (\ref{Stph}) is  just the Euler-Maclaurin
formula of Theorem \ref{LEMLF} and the second equation follows
from obvious estimates on the polynomial behaviour of the
integrals.

For $t=0$, as both sides of (\ref{Stph}) are quasipolynomials,
they take the same value. We may also deduce the equality for
$t=0$ from Corollary \ref{sommeconesausommet}. Indeed, for $t=0$,
the left hand side is $h(0)$ and the faces of dimension $>0$ give
a zero contribution to the right hand side. The equality becomes
$$
h(0)= \sum_{v\in\CV(\p)}(D(\dir(\t(\p,v)))\cdot h)(0),
$$
and follows immediately from Corollary \ref{sommeconesausommet}. $\Box$
\end{proof}

 For instance, let $h$ be a monomial $h(x)= x^m$ with
$m\in\N^d$. The coefficient of the highest degree term $t^{d+
\vert m\vert }$ is the integral $\int_{\p} x^m dx$.  As the
operator $D(\p,\p,t)$ is equal to $1$, the face $\p$  does not
contribute to the coefficients of lower degree. The
 coefficient of  $t^{d+ \vert m\vert -1}$ involves only the faces of
codimension $1$,  the
 coefficient of  $t^{d+ \vert m\vert -2}$ involves only the faces of
codimension $1$ and $2$, etc.

When we apply the last proposition to the function $h(x)=1$, we
obtain a formula for  the number of integral points in $t\p$.
\begin{corollary}\label{period}
(a) The Ehrhart quasipolynomial of the polytope $\p$ is given by
$$
\card(t\p\cap\lattice)=\sum_{\f\in\CF(\p)}
\nu_0(\p,\f,t)\vol(\f)t^{\dim\f},
$$
hence
$$
E_k(\p,t)= \sum _{\f, \dim \f =k} \nu_0(\p,\f,t)\vol(\f).
$$
(b)  The rational number $\nu_0(\p,\f,t)$ depends only on the class
modulo lattice translations of the transverse cone $t\t(\p,\f)$.
Therefore, it is a periodic function of $t$ with period at most
equal to $q_{\f}$, the smallest integer such that $q_{\f}<\f>$
contains integral points. In particular, if the affine span $<\f>$
contains integral points for every $k$-dimensional face $\f$  of
$\p$ , then the Ehrhart coefficient $E_k(\p,t)$  does not depend on
$t$.

(c) When $\dim V$ is fixed, there exists a polynomial time
algorithm which computes  $\nu_0(\p,\f,t)$.
\end{corollary}
\begin{proof}
The last statement in {\it (b)}  is due to Stanley \cite{Stanley};
it also follows immediately from Remark \ref{integralface}. The
computability follows from the computability of the functions
$\mu(\a)$. $\Box$  \end{proof}

Barvinok \cite{bb} proved recently  that, given an integer $m$,
there exists a polynomial time algorithm which computes the $m$
highest coefficients of the Ehrhart quasipolynomial of any
rational simplex in $\R^d$, \emph{when the dimension $d$ is
considered as an input}. We hope that our construction of the
functions $\mu$ leads to another polynomial time algorithm which
would compute the $m$ highest coefficients of the Ehrhart
quasipolynomial for any simplex in $\R^d$ and any polynomial
$h(x)$, when the dimension $d$ and the degree of the polynomial
are  considered as input.
\section{Computations in dimension 2}
In this section,  $V=\R^2$ with $\lattice=\Z^2$. Let $v_1$ and
$v_2$ be primitive integral vectors.  Let $s\in V_\Q$. We are
going to compute  $\mu(\a)$ for the affine cone
$$
\a=s+\R_+ v_1+\R_+ v_2.
$$
We will use the following notations:

\noindent For $t\in \R$,  we denote  the smallest
    integer greater or equal than $t$ by ${\rm ceil}(t)$ and we define $[[t]] \in [0,1[ $
by  $[[t]]= {\rm ceil}(t)-t$.

\noindent We denote by $B$  the function $B(y,t)=\frac{e^{[[t]]
y}}{1-e^y}+\frac{1}{y}$. Recall that  the function $\mu(\a)$ for a
one-dimensional cone (half-line) $\a$ is expressed in terms of $B$.

\noindent
$s=s_1 v_1 + s_2 v_2$ with $s_i\in\Q$.\\
$C_i= \frac{Q(v_1,v_2)}{Q(v_i,v_i)}$, for $i=1,2$.\\
$q= \det(v_1,v_2)$. We assume $q>0$.\\
$w\in \Z^2$ is a vector such that $\det(v_1,w)=1$,\\
$p=  \det(v_2,w)$. Thus $p$ and $q$ are coprime integers. \\
$r= (qs_1+[[qs_1]]) +p  (qs_2+[[qs_2]])$.\\
$\zeta$ is a primitive $q$-th root of $1$.

\noindent We observe that the lattice in  $ V/\R v_1$ is generated
by $\bar w= \frac{1}{q}\bar v_2$ where $\bar w$, $\bar v_2$ denotes
the image of $w$, $v_2$ in $ V/\R v_1$. Thus the transverse cone
$\t(\a,\f_1)\subset V/\R v_1$ is  given by
$$\t(\a,\f_1)= (qs_2+\R_+)\frac{1}{q}\bar v_2.$$

\begin{proposition}
For $\xi\in V^*$, let $y_i= \la \xi,v_i\ra$, for $i=1,2$. Then the
function $\mu(\a)(\xi)$ is given by
\begin{eqnarray}\label{munonunimodulaire}
\mu(\a)(\xi)=\hspace{8cm}
\end{eqnarray}
\begin{eqnarray*}
&&\frac{1}{q}
e^{[[qs_1]]\frac{y_1}{q}}e^{[[qs_2]]\frac{y_2}{q}}\left(
\frac{1}{(1-e^{\frac{y_1}{q}})(1-e^{\frac{y_2}{q}})}
 + \sum_{k=1}^{q-1} \frac{\zeta^{kr}}
{(1-\zeta^{k}e^{\frac{y_1}{q}}) (1-
  \zeta^{kp}e^{\frac{y_2}{q}})} \right)\\
 && + \frac{1}{y_1} B(\frac{y_2-C_1y_1}{q},[[q s_2]])
+\frac{1}{y_2 }B(\frac{y_1- C_2 y_2}{q},[[q s_1]])- \frac{ q}{y_1
y_2}.
\end{eqnarray*}
Its value at $\xi=0$ is equal to
\begin{eqnarray}\label{mu0nonunimodulaire}
 \mu(\a)(0)&= &\frac{1}{q}( (\half-[[qs_1]])(\half-[[qs_2]]))\\
\nonumber  +\frac{Q(v_1,v_2)}{Q(v_1,v_1)}
(\frac{1}{12}-\half[[qs_2]]& +&\half[[qs_2]]^2)
+\frac{Q(v_1,v_2)}{Q(v_2,v_2)}
(\frac{1}{12}-\half[[qs_1]]+\half[[qs_1]]^2) \\
\nonumber &+ &\frac{1}{q}  \sum_{k=1}^{q-1} \frac{\zeta^{kr}}{(1-\zeta^{k})
  (1-  \zeta^{kp})}.
\end{eqnarray}

When $\a$ is unimodular, $\mu(\a)(\xi)$ is given by
\begin{eqnarray}\label{unimodulaire}
\mu(a)(\xi)&= &\frac{\exp([[s_1]]y_1+[[s_2]]y_2)} {
(1-e^{y_1})(1-e^{y_2})} \\
\nonumber &+ &\frac{1}{y_1} B(y_2-C_1y_1,[[s_2]]) + \frac{1}{y_2 }
B(y_1- C_2 y_2,[[s_1]])- \frac{1}{y_1 y_2}
\end{eqnarray}
and its constant term $\mu(\a)(0)$ is given by
\begin{eqnarray}
&&\mu(\a)(0)=
(\half-[[s_1]])(\half-[[s_2]])\\
\nonumber &+ &\frac{Q(v_1,v_2)}{Q(v_1,v_1)}
(\frac{1}{12}-\half[[s_2]]+\half[[s_2]]^2) +\frac{Q(v_1,v_2)}{Q(v_2,v_2)}
(\frac{1}{12}-\half[[s_1]]+\half[[s_1]]^2).
\end{eqnarray}
\end{proposition}
\begin{remark}
In actual Maple computations, we use only the unimodular case
(\ref{unimodulaire}) which is computable in polynomial time at any
given order.  Thanks to the valuation property, for a non
unimodular cone $\a$, we compute $\mu(\a)$ by  performing first a
signed decomposition of $\a$ into unimodular cones, similar to
Barvinok's decomposition.  As a result, by our local
Euler-Maclaurin formula, we have fast algorithms which compute,
for a polygon $\p\subset\R^2$ and a monomial  $h(x)=x_1^{m_1}
x_2^{m_2}$, the sum of values at integral points $S(\p,h)=
\sum_{x\in \p\cap\lattice}h(x)$ and the coefficients of the
corresponding Ehrhart quasipolynomial.
\end{remark}
\begin{proof}
We use the defining relation of Proposition \ref{defdemu}. First,
we obtain  a summation formula for $S(\a)$ by using finite Fourier
transform as in \cite{bv1}. We observe that $\lattice \subset M=
\frac{1}{q}(\Z v_1+ \Z v_2)$. Let $\tilde \a\subset \a$  be the
cone
$$
\tilde \a= \tilde s +\R_+v_1 +\R_+v_2
$$
with vertex
$$
 \tilde s =\frac{1}{q} ({\rm ceil}(qs_1)v_1+{\rm ceil}(qs_2)v_2).
$$
Then $\tilde \a\cap M= \a\cap M$. As  $\lattice \subset  M$, we have also
 $\tilde \a\cap \lattice= \a\cap\lattice$.
Consider the dual lattice $M^*\subset \lattice^*$. Let $x\in M$.
We have
$$
\sum_{\gamma\in \lattice^* /M^*}e^{2i\pi \la \gamma,x\ra}= \begin{array}{l} 0
  \mbox{ if } x\notin \Z^2\\
 q  \mbox{ if } x\in \Z^2. \end{array}
$$
Therefore we have
$$
S(\tilde \a)(\xi)= \frac{1}{q} \sum_{\gamma\in \lattice^*
  /M^*}\sum_{x\in \tilde \a\cap M}e^{\la 2i\pi  \gamma+\xi,x\ra}.
$$
Since
$$
\tilde \a \cap M =  \tilde s + \Z_+ \frac{v_1}{q} +  \Z_+
\frac{v_2}{q},
$$
we obtain
\begin{equation}\label{Brionvergne}
S(\a)(\xi)= S(\tilde \a)(\xi)= e^{\la \xi,\tilde s \ra} \frac{1}{q}
\sum_{k=0}^{q-1}
\frac{e^{\la 2i\pi k \delta, \tilde s\ra}}
{ (1- e^{\la   2i\pi k\delta + \xi, \frac{v_1}{q}\ra} ) (1- e^{\la
  2i\pi k\delta + \xi, \frac{v_2}{q}\ra})}.
\end{equation}
In this formula, $\delta$ is a generator of the group
$\Lambda^*/M^* \cong \Z/q\Z$. By using the basis $(v_1,w)$ of
$\Z^2$, we obtain:
$$
\la\delta,v_1\ra =1 ,\;\;  \la\delta,v_2\ra =p.
$$
Let $(v_1^*, v_2^*)$ be the dual basis of $(v_1, v_2)$.
The orthogonal projection of $\xi $ on
$(V/\R v_1)^* = \R v_2^*$ is equal to
$(- C_1y_1+y_2)v_2^* $,  with
$$
C_1=- \frac{Q(v_1^*,v_2^*)}{Q(v_2^*,v_2^*)}
= \frac{Q(v_1,v_2)}{Q(v_1,v_1)}.
$$
Then  the computation in dimension one (\ref{mudim1}) gives
$$
\mu(\t(\a,\f_1))(\xi)= B(\frac{- C_1y_1+y_2}{q},[[qs_2]])
$$
and similarly
$$
\mu(\t(\a,\f_2))( \xi)= B(\frac{y_1 - C_2 y_2}{q},[[qs_1]])
$$
with
$$
C_2= \frac{Q(v_1,v_2)}{Q(v_2,v_2)}.
$$
We have
$$
I(\f_i)(\xi)= e^{y_1s_1+y_2s_2}(\frac{-1}{y_i}),
$$
$$
I(\a)(\xi)=  e^{y_1s_1+y_2s_2}\frac{q}{y_1y_2}.
$$
Therefore, by  (\ref{defdemu}), we have
 \begin{eqnarray}\label{mudim2}
 & &\mu(\a)(\xi) =
 e^{-(y_1s_1+y_2s_2)}S(\a)(\xi)\\
\nonumber & + &\frac{1}{y_1} B(\frac{-C_1y_1+y_2}{q},[[q s_2]])
+\frac{1}{y_2 }B(\frac{y_1- C_2 y_2}{q},[[q s_1]])\\
\nonumber & - & \frac{q}{y_1 y_2}.
\end{eqnarray}
In Formula (\ref{mudim2}), we replace $S(\a)$ with the right hand
side of (\ref{Brionvergne}), taking in account the equality ${\rm
ceil}(qs_i) -qs_i= [[qs_i]]$.
 This gives
(\ref{munonunimodulaire}). $\Box$
\end{proof}

 If $\a$ is not
unimodular, then $\mu(\a)(\xi) $  involves the ``extended''
Fourier-Dedekind sum
$$
 \frac{1}{q}\sum_{k=1}^{q-1} \frac{\zeta^{kr}}
{(1-  \zeta^{k}e^{\frac{y_1}{q}}) (1-
  \zeta^{kp}e^{\frac{y_2}{q}})},
$$
and  $\mu(\a)(0) $ involves the Fourier-Dedekind sum
$$
D(q,1,p,r) = \frac{1}{q}  \sum_{k=1}^{q-1} \frac{\zeta^{kr}}{(1-\zeta^{k})
  (1-  \zeta^{kp})}.
$$
One has (see for instance \cite{BDR})
$$
D(q,1,p,r)= \sum_{k=0}^{q-1}((-\frac{kp+r}{q}))((\frac{k}{q}))-\frac{1}{4q},
$$
where the `''sawtooth'' function $((a))$ is defined by
$$
((a))= a - [a]-\half.
$$
The Dedekind sum $D(q,1,p,r)$  can also be computed in polynomial
time by means of reciprocity relations (see for example
\cite{BDR}), but here we do not use this fact.

\begin{figure}[!h]
\begin{center}
\includegraphics{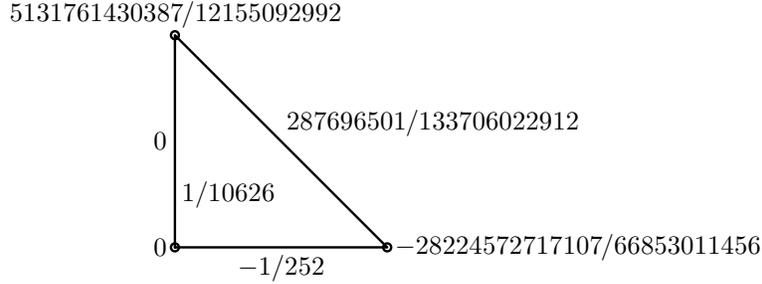}
\caption{Le savant Cosinus}\label{cosinus}
\end{center}
\end{figure}
\begin{example} (Figure\ref{cosinus})
We compute the right hand side of Euler-Maclaurin formula in the case
  of  the  ``dull triangle'' with vertices
  $(0,0)$,$(1,0)$,$(0,1)$, and the polynomial $h(x)=x_1^{20} x_2$.
As expected, the contributions of the various faces of $\p$  add up to
  $0$\footnote{This computation delighted us, and it would have delighted
Dr. Pancrace Eus\`ebe Z\'ephyrin Brioch\'e alias "Dr. Cosinus"
\cite{christophe}}.

\noindent Contribution of vertices: $0$,
$-{\frac{28224572717107}{66853011456}}$, ${\frac
{5131761430387}{12155092992}} $.

\noindent Contribution of edges:  $-{\frac {1}{252}}$, ${\frac
{287696501}{133706022912}}$, $0$.

\noindent Integral over triangle:  ${\frac {1}{10626}}$.
\end{example}
\begin{figure}[!h]
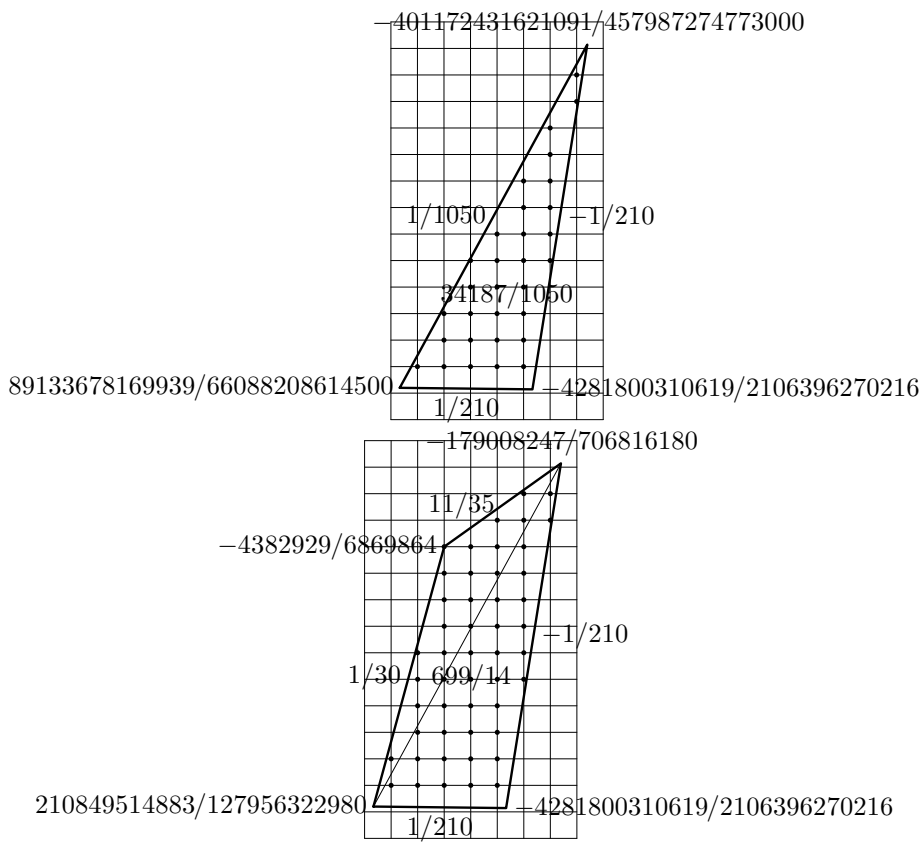

\begin{center}
\includegraphics{triangle357.01.eps}
\includegraphics{quadrilatere357.01.eps}
\caption{Triangle with vertices $(\frac{1}{3},\frac{1}{5})$,
$(\frac{16}{3},\frac{1}{7})$, $(\frac{37}{5},\frac{92}{7})$ and
quadrangle with extra vertex (3,10)}\label{nodelzant}
\end{center}
\end{figure}

\begin{example}\label{357}
Triangle with vertices
$s_1=(\frac{1}{3},\frac{1}{5})$, $s_2=(\frac{16}{3},\frac{1}{7})$,
$s_3=(\frac{37}{5},\frac{92}{7})$.

\noindent Number of integral points: $31$

\noindent Contribution of vertices: ${\frac
{89133678169939}{66088208614500}}$, $ -{\frac
{4281800310619}{2106396270216}}$, $-{\frac
{401172431621091}{457987274773000}}$

\noindent Contribution of edges :${\frac {1}{210}}, -{\frac
{1}{210}}, {\frac {1}{1050}}$

\noindent Area of  triangle: ${\frac {34187}{1050}}$
\end{example}
\begin{example}\label{quadrilatere357}
Quadrangle with vertices $(\frac{1}{3},\frac{1}{5})$,
$(\frac{16}{3},\frac{1}{7})$, $(\frac{37}{5},\frac{92}{7})$,
$(3,10)$.

\noindent Number of integral points: $49$

\noindent Contribution of vertices: ${\frac
{210849514883}{127956322980}}, -{\frac
{4281800310619}{2106396270216}}, -{\frac {179008247}{706816180}},
-{\frac {4382929}{6869864}}$

\noindent Contribution of edges : ${\frac {1}{210}}, -{\frac
{1}{210}}, {\frac {11}{35}}, {\frac {1}{30}}$

\noindent Area: ${\frac {699}{14}}$
\end{example}
Remark that, as expected, the contributions of the bottom right
vertex $(\frac{16}{3},\frac{1}{7})$ in the triangle or the
trapezoid of Figure \ref{nodelzant} are the same, as this vertex
have the same tangent cone in both polygons.

\begin{example}\label{ex357}
We compute the  Ehrhart quasipolynomial $E_2 t^2 + E_1(t)t +
E_0(t)$ for the number of integral points  of the triangle of
Example \ref{357}. The highest coefficient is the area of the
triangle, $E_2 = {\frac {34187}{1050}}$. The coefficient $E_1(t)$
is the sum of the contributions of the edges. The coefficient
$E_0(t)$ is the sum of the contributions of the vertices.

On this example, we can observe the periods of the contributions
of the edges and vertices to the Ehrhart coefficients (Corollary
\ref{period}). The period of  a vertex contribution is equal to
the \emph{lcm} of the denominators of its coordinates. For an edge
starting from a vertex $(a_1,a_2)$ and parallel to the primitive
vector $(v_1,v_2)$, the period is the least integer $q$ such that
$q(a_1v_2-a_2v_1)$ is an integer.

\medskip
\noindent Contribution of edges  (periods $3$, $5$, $7$
respectively):

\noindent $-{\frac {1}{105}}\,{\it mod} \left( t,3 \right) +{\frac
{1}{70}}, \;\;\; -{\frac {1}{105}}\,{\it mod} \left( 4\,t,7 \right)
+1/30, \;\;\; -{\frac {1}{525}}\,{\it mod} \left( 2\,t,5 \right)
+{\frac {1}{210}}$.

\medskip

\noindent Contribution of vertex $s_1=(\frac{1}{3},\frac{1}{5})$
(period $15$):
$$
\begin{array}{l}
{\frac {1}{75}}\,{\it mod} \left( t,5 \right) {\it mod} \left( t,15
 \right) -\frac{1}{45}\,{\it mod} \left( t,3 \right) {\it mod} \left( t,15
 \right) -\frac{1}{15}\,{\it mod} \left( t,5 \right) {\it mod} \left( 2\,t,3
 \right) \\
+{\frac {1765457034769}{293725371620}}-{\frac {44}{5}}\,{\it mod}
\left( t,5 \right) -{\frac {747989}{11987225}}\,{\it mod}
 \left( 2\,t,5 \right)
+\frac{1}{30}\,{\it mod} \left( 7\,t,15 \right) \\+\frac{3}{5}\,
{\it mod} \left( 8\,t,15 \right) -\frac{1}{15}\,{\it mod} \left(
t,15
 \right) +{\frac {44}{25}}\, \left( {\it mod} \left( t,5 \right)
 \right) ^{2}\\
-{\frac {2}{45}}\, \left( {\it mod} \left( 8\,t,15
 \right)  \right) ^{2}
+{\frac {1}{225}}\, \left( {\it mod} \left( t, 15 \right)  \right)
^{2}+\frac{1}{3}\,{\it mod} \left( 2\,t,3 \right)\\
 +{ \frac
{15227}{183774}}\,{\it mod} \left( t,3 \right) -\frac{1}{9}\, \left(
{ \it fmod} \left( 2\,t,3 \right)  \right) ^{2}+{\frac
{2567}{91887}}\,
 \left( {\it mod} \left( t,3 \right)  \right) ^{2}\\
-{\frac {901467}{ 119872250}}\, \left( {\it mod} \left( 2\,t,5
\right)  \right) ^{2}-{ \frac {1}{150}}\, \left( {\it mod} \left(
7\,t,15 \right) \right) ^{ 2}+{\frac {1}{75}}\,{\it mod} \left(
2\,t,5 \right) {\it mod}
 \left( 7\,t,15 \right) \\
+\frac{1}{45}\,{\it mod} \left( 8\,t,15 \right) {\it mod} \left(
2\,t,3 \right)
\end{array}
$$
\medskip

\noindent Contribution of vertex  $s_2=(\frac{16}{3},\frac{1}{7})$
(period $21$):
$$
\begin{array}{l}
-{\frac {32132693735}{4776408776}}-{\frac {13}{14}}\, \left( {\it
mod } \left( t,7 \right)  \right) ^{2}+\frac{13}{2}\,{\it mod}
\left( t,7
 \right)
+{\frac {1}{63}}\,{\it mod} \left( t,3 \right) {\it mod}
 \left( 2\,t,21 \right)\\
 -{\frac {15227}{183774}}\,{\it mod} \left( t,
3 \right)
 -\frac{1}{21}\,{\it mod} \left( 2\,t,21 \right)
-{\frac {2567}{ 91887}}\, \left( {\it mod} \left( t,3 \right)
\right) ^{2}+\frac{1}{9}\,
 \left( {\it mod} \left( 2\,t,3 \right)  \right) ^{2}\\
-\frac{2}{3}\,{\it mod}
 \left( 2\,t,3 \right)
+{\frac {1}{63}}\,{\it mod} \left( 2\,t,21
 \right) {\it mod} \left( 2\,t,3 \right)\\
+\frac{1}{21}\,{\it mod} \left( 2\, t,3 \right) {\it mod} \left(
4\,t,7 \right)   -{\frac {30189}{545804}} \,{\it mod} \left( 4\,t,7
\right)\\
 -{\frac {8797}{3820628}}\, \left( { \it mod} \left(
4\,t,7 \right) \right)^{2}
\end{array}
$$
\medskip

\noindent Contribution of vertex $s_3=(\frac{37}{5},\frac{92}{7})$
(period $35$):
$$
\begin{array}{l}
{\frac {9}{1225}}\, ( {\it mod} ( 16\,t,35 )  ) ^{2} +{\frac
{3}{1225}}\, ( {\it mod} ( 23\,t,35  )  ) ^{2} +{\frac
{641856910509}{373867163080}}\\
 +{\frac { 30189}{545804}}\,{\it mod}
( 4\,t,7 )
 +{\frac {1}{1225}}\,
 ( {\it mod} ( 34\,t,35 )  ) ^{2}\\
+1/10\,{\it mod} ( 3\,t,5 ) -1/35\,{\it mod} ( 34\,t,35 )
 -{\frac
{1}{2450}}\, ( {\it mod} ( 18\,t,35 ) ) ^{2}\\
 +{\frac {2}{1225}}\,
( {\it mod} ( 9\,t,35
 )  ) ^{2}
+{\frac {901467}{119872250}}\, ( {\it mod}  ( 2\,t,5 )  ) ^{2}\\
+{\frac {1}{70}}\,{\it mod}  ( 18\,t,35 ) +{\frac {1}{50}}\, ( {\it
mod} ( 3 \,t,5 )  ) ^{2}
 +{\frac {1}{1225}}\,{\it mod} ( 34\,t ,35
) {\it mod} ( 18\,t,35
)\\
 +{\frac {1}{1225}}\,{ \it fmod} ( 23\,t,35 ) {\it mod} ( 16\,t,35
)\\
 + {\frac {1}{175}}\,{\it mod} ( 3\,t,5 ) {\it mod} ( 23\,t,35 )
 -{\frac {1}{1225}}\,{\it mod} ( 34\,t,35
 ) {\it mod} ( 16\,t,35 )\\
-{\frac {1}{175}}\,{\it mod} ( 2\,t,5 ) {\it mod} ( 16\,t,35 )\\
-1/35\, {\it mod} ( 3\,t,7 ) {\it mod} ( 3\,t,5 )\\
 -1/ 35\,{\it mod} ( 3\,t,5 ) {\it mod} ( 4\,t,7 )
 -{\frac {1}{245}}\,{\it mod} ( 9\,t,35 ) {\it mod}  ( 3\,t,7 )\\
 -{\frac {1}{175}}\,{\it mod} ( 3\,t,5
 ) {\it mod} ( 9\,t,35 )\\
 -{\frac {1}{1225}}\,{\it mod} ( 16\,t,35 ) {\it mod} ( 18\,t,35 )
 +{ \frac {747989}{11987225}}\,{\it mod} ( 2\,t,5 )\\
 +{\frac { 8797}{3820628}}\, ( {\it mod} ( 4\,t,7 )  ) ^{2 }
-1/35\,{\it mod} ( 9\,t,35 )
 -{\frac {4}{35}}\,{\it mod}( 23\,t,35 )\\
 -{\frac {8}{35}}\,{\it mod} ( 16\,t,35  )
+1/49\, ( {\it mod} ( 3\,t,7 )  ) ^{2}
\end{array}
$$
\end{example}
\begin{example}({\bf Computation time}).

We computed the full Ehrhart quasipolynomial corresponding to the
triangle of Example \ref{357} and  the polynomial $h(x_1,x_2)=
x_1^{k_1}x_2^{k_2}$,  with increasing exponents $k_1$ and $k_2$.
Allowing a computation time of about one hour, we reached $k_1=k_2=
24$. The result is of course too big to write  here.

The sum of values $x_1^{48}x_2^{48}$ at the integral points of the
triangle of Example \ref{357} dilated by the  factor $N=11^5$ took
about the same time. The result is the following number
$$
\begin{array}{l}
55969247458735493271268368615238071121335974262337882261418363621\\
89704055956429496253759473056373507451253522021344188115187647607\\
84555431172202923756940824265247663088847763429436570335188702325\\
06644969965841257822711805056447218921550669146263582661876630783\\
21357671611262065293901983868557252464459832189159990869820527095\\
53646871654914800005753059422066576204781923454823934475242960034\\
42199041253798398004263030681714027295470241663946228744550160085\\
43856624239377702107746492579014275563017167813144052693763385569\\
75239252588060279466314599314734680953729093269435217987689840619\\
0740089242444014302.
\end{array}
$$
\end{example}

As experiments showed,  our method for the computation of $\sum_{x
\in \p \cap \Z^d} h(x)$ is very efficient for this small dimension,
compared to other available softwares. Furthermore, as Example
\ref{ex357} shows, the Ehrhart polynomial is written as a sum of
canonical contributions of all faces, once the scalar product is
fixed. We will come back soon to the computational and complexity
aspects of this problem for higher dimensions.

\end{document}